\documentclass[12pt]{article}

\usepackage[colorlinks=true,
linkcolor=webgreen,
filecolor=webbrown,
citecolor=webgreen]{hyperref}
\definecolor{webgreen}{rgb}{0,.5,0}
\definecolor{webbrown}{rgb}{.6,0,0}

\usepackage{amsmath,amssymb,graphics,epsfig}

\makeatletter
\renewcommand{\thesection}{\@Roman{\c@section}}
 \renewcommand{\thesubsection}{\@Alph{\c@subsection}}
 \renewcommand{\thesubsubsection}{\thesubsection\@arabic{\c@subsubsection=
}}
 \def\@secnum{\thesection}
\makeatother

\title{{\sffamily
Asymmetric Multiple Description Lattice Vector Quantizers
}}
\author{ 
\begin{tabular}{cp{10mm}cp{10mm}c}
Suhas N. Diggavi, N. J. A. Sloane, and Vinay A. Vaishampayan 
\end{tabular}\\
AT\&T Shannon Laboratories,\\
180 Park Avenue, Bldg 103,\\
Florham Park NJ 07932, USA. \\
{\sffamily Tel: (973)-360-8492, FAX: (973)-360-8178.} \\
{\sffamily Email: \{suhas,njas,vinay\}@research.att.com}
}

\def\Reals{\mathop{\hbox{\mit I\kern-.2em R}}\nolimits}

\def\RA{\mathop{\hbox{\Rightarrow\kern-.2em /}}}

\def\Complexes{\mathop{\hbox{\mit C\kern-.44em
               \vrule depth 0ex height 1.4ex width .06em
               \kern.41em}}\nolimits}

\def\Zeals{\mathop{\hbox{\mit Z\kern-.29em Z}}\nolimits}

\def\Neals{\mathop{\hbox{\mit I\kern-.2em N}}\nolimits}

\def\K{\mathop{\hbox{\mit I\kern-.2em K}}\nolimits}

\def\F{\mathop{\hbox{\mit I\kern-.2em F}}\nolimits}

\def\sqr#1#2{{\vcenter{\vbox{\hrule height.#2pt
  \hbox{\vrule width.#2pt height#1pt \kern#1pt
    \vrule width.#2pt}
   \hrule height.2pt}}}}

\def \bX{{\bf X}}

\def \calE{{\cal E}}

\def \calH{{\cal H}}

\begin{document}

\thispagestyle{myheadings}
\maketitle
\setcounter{equation}{0}
\setcounter{figure}{0}

\pagenumbering{arabic}


%
%
\newcommand{\beq}{\begin{equation}}
\newcommand{\bear}{\begin{eqnarray}}
\newcommand{\eear}{\end{eqnarray}}
\newcommand{\bi}{\begin{itemize}}
\newcommand{\ei}{\end{itemize}}
\newcommand{\ben}{\begin{enumerate}}
\newcommand{\een}{\end{enumerate}}
\newcommand{\eeq}{\end{equation}}
\newcommand{\beql}[1]{\begin{equation}\label{#1}}
\newcommand{\eqn}[1]{(\ref{#1})}

\newtheorem{theorem}{Theorem}[section]
\newtheorem{corollary}{Corollary}[section]
\newtheorem{lemma}{Lemma}[section]
\newtheorem{property}{Property}
\newtheorem{claim}{Claim}[section]
\newtheorem{proposition}{Proposition}[section]
\newtheorem{definition}{Definition}[section]
\newtheorem{fact}{Fact}[section]

\newcommand{\half}{\raisebox{1ex}{\footnotesize 1}/\raisebox{-1ex}{\footnotesize 2}}
\newcommand{\lalalaU}{\Lambda_\cup = \langle \Lambda_1 , \Lambda_2 \rangle}
\newcommand{\lalalaA}{\Lambda_\cap = \Lambda_1 \cap \Lambda_2}
\newcommand{\abf}{\mbox{${\bf a }$} }
\newcommand{\cbf}{\mbox{${\bf c }$} }
\newcommand{\chbf}{\mbox{${\bf \hat{c} }$} }
\newcommand{\ebf}{\mbox{${\bf e }$} }
\newcommand{\fbf}{\mbox{${\bf f }$} }
\newcommand{\gbf}{\mbox{${\bf g }$} }
\newcommand{\hbf}{\mbox{${\bf h }$} }
\newcommand{\hhbf}{\mbox{${\bf \hat{h}}$} }
\newcommand{\nbf}{\mbox{${\bf n }$} }
\newcommand{\sbf}{\mbox{${\bf s }$} }
\newcommand{\wbf}{\mbox{${\bf w }$} }
\newcommand{\xbf}{\mbox{${\bf x }$} }
\newcommand{\ybf}{\mbox{${\bf y }$} }
\newcommand{\zbf}{\mbox{${\bf z }$} }

\newcommand{\Abf}{\mbox{${\bf A }$} }
\newcommand{\Bbf}{\mbox{${\bf B }$} }
\newcommand{\Dbf}{\mbox{${\bf D }$} }
\newcommand{\Ebf}{\mbox{${\bf E }$} }
\newcommand{\Fbf}{\mbox{${\bf F }$} }
\newcommand{\Gbf}{\mbox{${\bf G }$} }
\newcommand{\Hbf}{\mbox{${\bf H }$} }
\newcommand{\Ibf}{\mbox{${\bf I }$} }
\newcommand{\Kbf}{\mbox{${\bf K }$} }
\newcommand{\Mbf}{\mbox{${\bf M }$} }
\newcommand{\Nbf}{\mbox{${\bf N }$} }
\newcommand{\Pbf}{\mbox{${\bf P }$} }
\newcommand{\Qbf}{\mbox{${\bf Q }$} }
\newcommand{\Rbf}{\mbox{${\bf R }$} }
\newcommand{\Sbf}{\mbox{${\bf S }$} }
\newcommand{\Ubf}{\mbox{${\bf U }$} }
\newcommand{\Vbf}{\mbox{${\bf V }$} }
\newcommand{\Xbf}{\mbox{${\bf X }$} }
\newcommand{\Ybf}{\mbox{${\bf Y }$} }
\newcommand{\Zbf}{\mbox{${\bf Z }$} }
\newcommand{\ZZ}{\mathbb Z}
\newcommand{\sP}{\cal P}

\newcommand{\alphabf}{\mbox{${\boldmath \alpha }$} }
\newcommand{\betabf}{\mbox{${\boldmath \beta }$} }
\newcommand{\epsibf}{\mbox{${\boldmath \epsilon }$} }
\newcommand{\psibf}{\mbox{${\boldmath \psi }$} }
\newcommand{\xibf}{\mbox{${\boldmath \xi }$} }
\newcommand{\thetabf}{\mbox{${\boldmath \theta }$} }
\newcommand{\Thetabf}{\mbox{${\boldmath \Theta }$} }
\newcommand{\Lambdabf}{\mbox{${\boldmath \Lambda }$} }
\newcommand{\Zerobf}{\mbox{${\large{\bf 0 }}$} }
\newcommand{\La}{\Lambda}
\newcommand{\la}{\lambda}

\newcommand{\Expt}{\mbox{${\Bbb E}$} }
\newcommand{\Expth}{\mbox{$\hat{{\Bbb E}}$} }
\newcommand{\gauss}{\mbox{${\cal N}$} }
\newcommand{\Prob}{\mbox{${\mathbb P}$} }
\newcommand{\expect}[1]{\Expt \left[ #1 \right]}
\newcommand{\prob}[1]{\Prob \left\{ #1 \right\}}

\newcommand{\indr}[1]{{\bf 1}_{\left\{ #1 \right\}}}
\newcommand{\past}[1]{ {#1}^{i-1}}
\newcommand{\Qpast}{\past{Q}}
\newcommand{\Xpast}{\past{X}}
\newcommand{\Ypast}{\past{Y}}
\newcommand{\pnsub}[2]{p \left( #2 \right) } 
\newcommand{\psub}[2]{p_{#1} \left( #2 \right) } 
\newcommand{\pg}[2]{\pnsub{#1}{#2}}
\newcommand{\pyicqiyp}[3]{\pg{Y_i | Q_i, \Ypast}{#1 | #2, #3}} 
\newcommand{\pqicyp}[2]{\pg{Q_i | \Ypast}{#1 | #2}} 
\newcommand{\pqicqp}[2]{\pg{Q_i | \Qpast}{#1 | #2}} 
\newcommand{\pyicqixi}[3]{\pg{Y_i | Q_i, X_i}{#1 | #2, #3}}
\newcommand{\pqicypxixp}[4]{\pg{Q_i | \Ypast, X_i, \Xpast}{#1 | #2, #3, #4}}
\newcommand{\pyicqixixp}[4]{\pg{Y_i | Q_i, X_i, \Xpast}{#1 | #2, #3, #4}}
\newcommand{\px}[1]{\pg{X}{#1}} 

\newcommand{\Complex}{{\rm {\bf C}\mkern-9mu\rule{0.05em}{1.4ex}\mkern10mu}}
\newcommand{\be}{\begin{equation}}
\newcommand{\ee}{\end{equation}}
\newcommand{\benonum}{\begin{displaymath}}
\newcommand{\eenonum}{\end{displaymath}}
\def\bm#1{\mbox{\boldmath $#1$}}

\newenvironment{proof}{
	\noindent
	{\bf Proof:}
}{
	\hfill$\blacksquare$
}

\small
\subsection*{\centering Abstract}
{\em
We consider the design of asymmetric multiple description lattice quantizers
that cover the entire spectrum of the distortion profile, ranging from
symmetric or balanced to successively refinable. 
We present a solution to a labeling problem, which is an important part of
the construction, along with a general design procedure.
This procedure is illustrated using the $\ZZ^2$ lattice.
The asymptotic performance of the quantizer is analyzed in the high-rate case.
We also evaluate its rate-distortion performance and compare
it to known information theoretic bounds.
}
\normalsize


\section{Introduction}
\label{sec:intro}

A multiple description source encoder generates a set of binary
streams or descriptions of a source sequence, each with its own
rate constraint. The transmission medium may deliver some or all
of the descriptions to the decoder. The objective is to minimize
the distortion between the source sequence and the decoded sequence
when all the descriptions are available, while ensuring that the distortion
which results when only a subset of the descriptions are available remains
below a pre-specified value that depends on the subset.
If there are $D$ descriptions, the {\em distortion profile}
is a vector of length $2^D$ whose components give the distortion constraints for each subset of the descriptions.

In recent years, multiple description coders have received considerable
attention, driven by the interest in packet voice and
video communications
(see the bibliography).
Most of the work
(with the exception of \cite{FleEff:1})
has centered around the
successively refinable case and the balanced/symmetric case,
which are in a
sense two extremes of the distortion profile.
Successive refinement coders
find application in networks with a priority structure whereas
balanced codes are useful in networks that do not have such a structure, the
best example
at the present time being the Internet.

In this paper we propose a structured scheme that bridges the two cases,
in the sense that it
permits a fairly general distortion profile to be specified.
By allowing the individual descriptions to have different distortions,
the quantizer behavior can range from the balanced case
(where each description is equally important) to a strict hierarchy
(where the loss of some descriptions could make decoding impossible).
The new design is described in terms of a lattice vector quantizer, but the general
principle of asymmetric multiple description coding can be
extended to many other quantizers,
such as trellis coded quantizers, unstructured vector quantizers,
etc.
This could potentially allow us to incorporate channel (or network
route) reliability information into the transmission.
Also, it might be a useful way to allow for less intrinsic wastage of
network traffic as some descriptions could be given to the decoder
without necessarily waiting for the more important descriptions to arrive
(as in successive refinement).

For previous work on the information theoretic aspects of the multiple
description
problem see \cite{ECII,  EQCOI, OZAI, WWZI, ZBI}.
The problem of designing quantizers for the multiple
description problem has been considered in
\cite{FleEff:1, JafTar:2, MohRisLad:1, VAI2, VAI3, VaiSloSer}.
The work presented here extends that in
\cite{VaiSloSer}, which considered only the balanced/symmetric case.
Unlike the work in \cite{FleEff:1}, we do not use a training
approach; instead we use the geometry of the underlying
lattice to solve a labeling problem.
Other approaches to multiple description  coding based on overcomplete
expansions
are presented in \cite{BalDauVai:1, ChouMehWang:1, GoyKovVet:1}
and methods based on optimizing transforms and predictors are presented in
\cite{IVI,OWVR,VSI}.

The paper is organized as follows.
The source coding problem is formulated in
Section \ref{sec:formulation}, the design method is described in
Section \ref{sec:algo},
properties of the lattices and sublattices needed for
the construction are developed in Section \ref{sec:good}, 
a high rate analysis is
presented in Section \ref{sec:highrate}, and 
numerical results are presented in
Section~\ref{sec:results}.

\section{Preliminaries}
\label{sec:formulation}

\begin{figure}[h]
\centering
\input{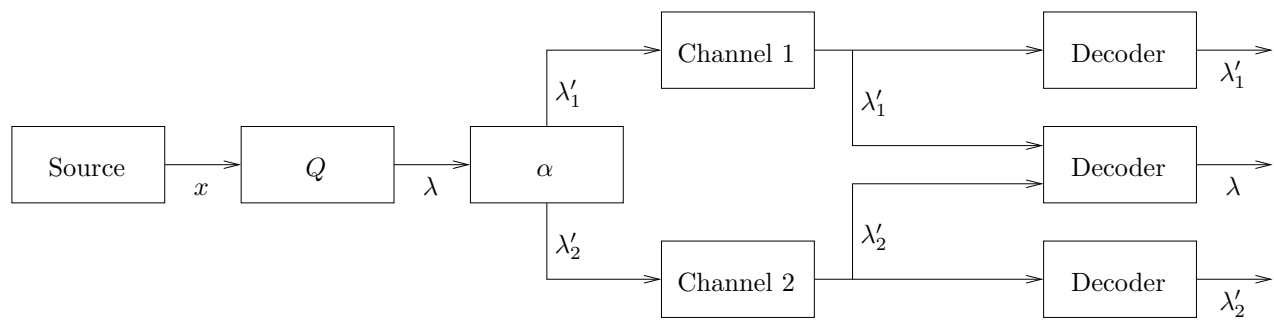}
\caption{Block diagram of a multiple description vector
quantizer.}
\label{fig-mdvq}
\end{figure}

A block diagram of a two-channel multiple description vector quantizer
(MDVQ) using a lattice codebook is shown in
Fig.~\ref{fig-mdvq}. An $L$-dimensional source vector $x$ is
first encoded as the closest vector
$\lambda$ in a lattice $\Lambda\subset \mathbb R ^L$. We will write
$\lambda=Q(x)$.
 Information about the selected code vector $\lambda$ is then
sent across the two channels, subject to rate constraints imposed
by the individual channels.
This is done through a labeling function $\alpha$.
At the decoder, if only channel 1 works, the received
information is used to select a vector $\lambda'_1$
from the channel 1 codebook. If only channel 2 works, the
information received over channel 2 is used to select
a code vector $\lambda'_2$ from the channel 2 codebook. 
If both channels work, it is assumed that 
enough information is available to recover $\lambda$.

We will assume that the channel 1 and channel 2 codebooks,
denoted by $\Lambda_1$ and $\Lambda_2$ respectively, are
sublattices\footnote{Strictly speaking, the codebooks are finite subsets of the
sublattices $\Lambda_1$ and $\Lambda_2$,
but we will ignore that distinction in this paper.}
of $\Lambda$.  
The index $[\Lambda: \Lambda_i ]$ is denoted by $N_i$, $i=1,2$.
$N_i$ is also called the {\em re-use index} 
of sublattice $\Lambda_i$.
We assume that each $\Lambda_i$ is {\em geometrically
similar} to $\Lambda$,
i.e.
that $\Lambda_i$ can be obtained from $\Lambda$ by applying a {\em similarity}
(a rotation, change of scale and possibly a reflection).
To simplify the analysis we will usually assume that the $\Lambda_i$ are
{\em strictly similar} to $\Lambda$, i.e. that reflections are not used.

\begin{property}\label{proii1}
Let $\Lambda$ be an $L$-dimensional lattice with generator matrix $G$ (the rows of $G$ span $\Lambda$).
A sublattice $\Lambda_1 \subseteq \Lambda$ is geometrically strictly similar to $\Lambda$ if and only if the following condition holds:
there is an invertible $L \times L$ matrix $U_1$ with integer entries,
a scalar $c_1$, and an orthogonal $L \times L$ matrix $K_1$ with determinant $1$ such that a generator matrix for $\Lambda_1$ can be written as
\begin{equation}\label{EqU1}
G_1 = U_1 G = c_1 G K_1 ~.
\end{equation}
If \eqn{EqU1} holds then the index of $\Lambda_1$ in $\Lambda$ is equal to
\begin{eqnarray}\label{EqU2}
N_1 & = & [ \Lambda : \Lambda_1 ] = \sqrt{\frac{\det \Lambda_1}{\det \Lambda}} =
\frac{\det G_1}{\det G} \nonumber \\
& = & \det U_1 = c_1^L ~.
\end{eqnarray}
Furthermore, $\Lambda_1$ has Gram matrix
\begin{equation}\label{EqU3}
A_1 = G_1 G_1^{tr} = U_1 G G^{tr} U_1^{tr} = U_1 AU_1^{tr} = c_1^2 A ~,
\end{equation}
where $A = G G^{tr}$ is a Gram matrix for $\Lambda$.
\end{property}

Even if the similarity is not strict, equations \eqn{EqU1}, \eqn{EqU2} and \eqn{EqU3} still hold but with $\det K_1 = -1$.

We will also usually assume that the sublattices $\Lambda_1$ and $\Lambda_2$ are {\em clean}
\cite{sim}, that is, no point of $\Lambda$ lies on the boundary of the Voronoi cells of $\Lambda_1$ or $\Lambda_2$.
Our algorithm still applies if this condition is not satisfied, but the book-keeping becomes more complicated.

Finally, we
require a sublattice 
$\Lambda_s$ of $\Lambda_1\bigcap \Lambda_2$
which is geometrically strictly similar
to $\Lambda$ and has index
$N_s=N_1N_2$ in $\Lambda$.
To reduce the complexity of the design we will also sometimes make use of
a sublattice
$\Lambda_{lcm}$ of $\Lambda_1\bigcap \Lambda_2$ 
which has index $N_{lcm}={\rm lcm}(N_1,N_2)$ in $\Lambda$
(such sublattices do not always exist -- see
Section \ref{sec:good}).

Since the information
sent over channel 1 is used to identify a
code vector $\lambda_1\in \Lambda_1$, and the information over channel 2 is
used to identify a code vector $\lambda_2\in \Lambda_2$,
we will assume that the labeling function $\alpha$
is a mapping from $\Lambda$ into
$\Lambda_1 \times \Lambda_2$ and that $(\lambda_1,\lambda_2)=\alpha(\lambda)$.
The component mappings are $\lambda_1=\alpha_1(\lambda)$ and
$\lambda_2=\alpha_2(\lambda)$.
In order to recover $\lambda$ when both channels
work, it is necessary that $\alpha$ be one-to-one.

Given $\Lambda$, $\Lambda_1$, $\Lambda_2$ and $\alpha$, there are three
distortions and two rates associated with the quantizer.
For a given source vector $x$ mapped to the triple $(\lambda,
\lambda_1,\lambda_2)$, the {\em two-channel
distortion} $d_0$ is given by $\|x-\lambda\|^2$, the side 
distortions $d_i$ by $\|x-\lambda_i\|^2, \,\,i=1,2$,
where $||\xbf||^2\stackrel{{\rm def}}{=}(1/L)\sum_{i=1}^Lx_i^2$
is the dimension-normalized Euclidean norm.
The corresponding average distortions are
denoted by $\bar{d}_0$,
$\bar{d}_1$ and $\bar{d}_2$.
(We will also refer to $\bar{d}_0$ as the {\em central} distortion.)
We assume that an entropy coder is used to transmit the
labeled
vectors at a rate arbitrarily close to the entropy, i.e.,
$R_i=\calH( \alpha_i(Q(\bX)) )/L$,
$i=1,2$, where
$\calH$ is
the binary entropy function.
The problem is to design the labeling function $\alpha$
so as to minimize $\bar{d}_0$ subject to $\bar{d}_1\leq D_1$ and
$\bar{d}_2\leq D_2$, for given rates
$(R_1,R_2)$ and
distortions $D_1$ and $D_2$.

We will assume that the source is memoryless with probability density function (pdf) $p$. The $L$-fold pdf
will be denoted by $p_L$ where $p_L((x_1,x_2,\ldots,x_L))=\prod_{i=1}^Lp(x_i)$. The 
differential entropies satisfy the relation $h(p_L)=Lh(p)$.

Given  a lattice $\Lambda$, a sublattice $\La'$ and a point $\lambda' \in \La'$,
we denote by  $V_{\Lambda : \La'} (\lambda' )$
the set of all points in $\Lambda$ that are
closer to $\lambda'$ than to any other point in $\La'$.
This set is the discrete Voronoi set of $\lambda'$ in $\Lambda$.
If $\La'$ is a clean sublattice of $\Lambda$ we do not need to worry
about ties when calculating $V_{\La : \Lambda'} (\lambda' )$.
The Voronoi cell $V_{\Lambda}(\lambda)$ of a point $\lambda \in \Lambda$ is
the set of all points
in $\mathbb R ^L$ that are
at least as close to $\lambda$ as to any other point of $\Lambda$.
Also $\calE(\lambda')=\alpha(V_{\Lambda: \La'}(\lambda'))$,
$\lambda' \in \La'$, will denote the set of all labels of the
points in $V_{\Lambda: \La'}(\lambda')$.

\subsection{\bf Distortion Computation}
The average two-channel distortion $\bar{d}_0$ is given by
\begin{equation}
\bar{d}_0=\sum_{\lambda \in
\Lambda}\int_{V_{\Lambda}(\lambda)}\|x-\lambda\|^2 p_L(x) dx.
\label{eqn-app-2}
\end{equation}
Since the codebook of the quantizer is a lattice, all the Voronoi
sets in the above summation are congruent. Furthermore,
upon assuming that each Voronoi cell is small and letting $\nu$ denote the $L$-dimensional
volume of a Voronoi cell, we obtain the two-channel distortion
\begin{equation}
\bar{d}_0=\frac{\int_{V_{\Lambda}(0)}\|x\|^2dx}{\nu} 
= G(\Lambda)\nu^{2/L},
\end{equation}
where the normalized second moment $G(\Lambda)$ is defined by
(\cite{SPLAG})
\begin{equation}
G(\Lambda)=\frac{\int_{V_{\Lambda}(0)}\|x\|^2 dx}{\nu^{1+2/L}}~.
\label{eqn-cen-dis}
\end{equation}

When only description $i$ is available, for $i=1,2$, the distortion is given by
\begin{equation}
\bar{d}_i =
 \bar{d}_0 + \sum_{\lambda \in \Lambda} \|\lambda-\alpha_i(\lambda) \|^2 P(\lambda),
\label{eqn-app-3}
\end{equation}
where $P(\lambda)$ is the probability of lattice point $\lambda$, and we have
assumed that $\lambda$ is the {\em centroid} of its Voronoi cell. This is true for the
uniform density. For nonuniform densities, there is an 
error term which goes to zero with the size of the Voronoi cell.
The first term in (\ref{eqn-app-3}) is the two-channel distortion and the second term is the excess
distortion
which is incurred when only description $i$ is available.
Note that, for a given $\Lambda$, only the excess distortion term
is affected
by the labeling $\alpha$.

At this point we impose a constraint on the labeling function that allows us to
reduce the problem to that of labeling a finite number of points.
We assume that the labeling function has the property that
$\alpha(\lambda+\lambda_s)=\alpha(\lambda)+\lambda_s$, for all $\lambda_s \in \Lambda_s$.
This leads to the following simplification:
\begin{equation}
\label{eq:DiComp}
\bar{d}_i=\bar{d}_0+(1/N_s)\sum_{\lambda \in
V_{\Lambda:\Lambda_s}(0)}\|\lambda-\alpha_i(\lambda)\|^2,
\end{equation} 
where we have assumed that $P(\lambda)$ is approximately constant
over a Voronoi cell of the sublattice $\Lambda_s$,
but may vary from one Voronoi cell to another.

\subsection{\bf Rate Computation}

Let $R_0$ bits/sample be
the rate required to address the
two-channel codebook for a single channel system\footnote{This quantity is
useful for evaluating the 
two-channel distortion as well as for evaluating the rate overhead associated
with the multiple description
system.}.
We first derive an expression for $R_0$ and then determine the rates
$R_1$ and $R_2$.
We use the fact that each quantizer bin
has identical volume $\nu$ and that $p_L(x)$ is approximately piecewise constant over
each Voronoi cell of $\Lambda_1$ and $\Lambda_2$.
This assumption is valid in the limit as the Voronoi cell become
small and is standard
in asymptotic quantization theory.

The rate $R_0=\calH(Q(X))$ is given by
\begin{eqnarray}
R_0 & = & -(1/L)\sum_{\lambda}\int_{V_{\Lambda}(\lambda)}p_L(x)dx \log_2
\int_{V_{\Lambda}(\lambda)} p_L(x) dx \nonumber \\
  & \approx & -(1/L) \sum_{\lambda}\int_{V_{\Lambda}(\lambda)}p_L(x)dx \log_2
p_L(\lambda) \nu \nonumber \\
  & \approx & h(p)-(1/L)log_2(\nu).
\label{eqn-rate-0}
\end{eqnarray}

It can be shown that the rate for description $i$ is given by
\begin{equation}
\label{eq:RiComp}
R_i=R_0-(1/L)\log_2(N_i),\quad i=1,2 \,.
\end{equation}
A single channel system would have used $R_0$ bits/sample. Instead a multiple
description system uses a
total of $R_1+R_2=2R_0-(1/L)\log_2(N_1N_2)$ bits/sample, and so
the rate overhead is $R_0-(1/L)\log_2(N_1N_2)$.

\section{Construction of the Labeling Function}
\label{sec:algo}

Suppose $\La$ is an $L$-dimensional
lattice with a pair of geometrically strictly similar, clean sublattices
$\La_1$ and $\La_2$,
and let $\La_s$ (the product sublattice)
be a geometrically strictly similar, clean sublattice
of both $\La_1$ and $\La_2$, with indices 
$[\La : \La_1] = N_1$, 
$[\La : \La_2] = N_2$ and $[\La : \La_s] = N_1 N_2$. 

In order to construct a labeling function we first identify
$\calE$, the subset of points of
$\Lambda_1\times \Lambda_2$ that will be used to label the points of $\Lambda$.
Next, a one-to-one correspondence will be established
between
$V_{\La:\La_s}(0)$ and a proper subset of $\calE$ so as
to minimize
an appropriate objective function, while ensuring that the labeling can
be extended
uniquely to the entire lattice.
To this end we first start by formulating a cost criterion
that will be used in the design.

\subsection{Cost criterion}

The multiple descriptions problem may be formulated \cite{ECII} as
a problem of
minimizing the central distortion subject to constraints on the
side distortion. 
The associated Lagrangian cost criterion is given by
\begin{eqnarray}
\label{eq:Lag}
J &=& \bar{d}_0 + \sum_{i=1}^2 \gamma_i \bar{d}_i
\\ \nonumber
&=& (\gamma_1+\gamma_2+1)\bar{d}_0 + \sum_{i=1}^2\gamma_i\sum_{\lambda \in \Lambda} 
\|\lambda-\alpha_i(\lambda) \|^2 P(\lambda)
\\ \nonumber
&=& (\gamma_1+\gamma_2+1)\bar{d}_0 + \sum_{\lambda \in \Lambda} P(\lambda)
\sum_{i=1}^2 \gamma_i\|\lambda-\alpha_i(\lambda) \|^2.
\end{eqnarray}
where $\gamma_1$, $\gamma_2$ are Lagrange multipliers.

The central distortion $\bar{d}_0$ is determined by
the lattice $\Lambda$.
If we assume that $P(\lambda)$ is approximately constant
over the Voronoi cell of $\Lambda_s$, we can rewrite the
cost criterion in terms of the cost over a Voronoi cell of $\Lambda_s$.
Then the design problem reduces to finding
a labeling scheme $\alpha(\lambda)$ which minimizes
\begin{equation}
\label{eq:opt_cr}
\frac{1}{N_s}\sum_{\lambda \in V_{\La:\La_s}(0)}
\left [ \gamma_1\|\lambda-\alpha_1(\lambda) \|^2 + 
\gamma_2\|\lambda-\alpha_2(\lambda) \|^2 \right ].
\end{equation}
After some algebra, the expression inside the summation can be rewritten as
\begin{equation}\label{eq:guide}
\frac{\gamma_1\gamma_2}{\gamma_1+\gamma_2}
\|\alpha_2(\lambda)-\alpha_1(\lambda)\|^2 + 
(\gamma_1+\gamma_2)\|\lambda-\frac{\gamma_1\alpha_1(\lambda)+\gamma_2
\alpha_2(\lambda)}{\gamma_1+\gamma_2}\|^2.
\end{equation}
The values of $\gamma_1$ and $\gamma_2$ determine the relative values of the two side-distortions $\bar{d}_1$ and $\bar{d}_2$.

Therefore our design principle is (informally)
for a given pair $\gamma_1$ and $\gamma_2$, to find a
labeling function $\alpha(\lambda)$ such that the sublattice
points $\alpha_1(\lambda)\in\Lambda_1,
\alpha_2(\lambda)\in\Lambda_2$ are not very far apart and
the lattice point $\lambda \in \La$ that is being labeled is not very far from the 
weighted mean
(the second term of (\ref{eq:guide})) of these two sublattice points.
This general guiding principle leads to our lattice design.
We will first describe the basic quantizer design and then
illustrate it using the lattice
$\ZZ^2$.

\subsection{Lattice Quantizer}
\label{subsec:lattice}

The quantizer construction is based on the following steps.

\begin{enumerate}

\item We are given an $L$-dimensional lattice $\Lambda$, rates $R_1,R_2$ and
distortions $D_1,D_2$.
These determine the indices $N_1,N_2$ using (\ref{eq:RiComp}), and
we attempt to find
(strictly similar, clean) sublattices
$\Lambda_1,\Lambda_2$ with these indices, together with
a product sublattice $\Lambda_s$.
We also choose appropriate values for
the weights $\gamma_1$ and $\gamma_2$.
For example, a successively refineable quantizer
corresponds to choosing $\gamma_1 =1$ and $\gamma_2 =0$.
For the balanced case we take $\gamma_1 = \gamma_2$.
By appropriately choosing $N_1,N_2,\gamma_1,\gamma_2$, 
one can achieve different levels of asymmetry in rate and distortion.

\item
We find the discrete Voronoi set\footnote{We usually
omit the word ``discrete'' when referring to this set.}
$V_0=V_{\La:\La_s}(0)$ 
for the sublattice $\Lambda_s$.
This is the fundamental set of points that we will label.
The labeling is then extended to the full lattice using the shift invariance
property (see Section \ref{sec:formulation}). 
We also find the sets
\begin{eqnarray}\label{EqP1}
{\cal P}_1 & = & V_{\La_1:\La_s}(0) ~=~
V_{\La:\La_s} (0) \cap \Lambda_1 ~, \\
\label{EqP2}
{\cal P}_2 & = & V_{\La_2:\La_s}(0) ~=~
V_{\La:\La_s} (0) \cap \Lambda_2 ~,
\end{eqnarray}
which are the
points of $\Lambda_1$ and $\Lambda_2$ belonging to the discrete Voronoi set.

\item
We determine the set 
\begin{equation}\label{EqP3}
{\cal L}_1(\lambda_1)=\{\lambda_2\in\Lambda_2:
\lambda_2\in V_0+\lambda_1\}
\end{equation}
for all $\lambda_1\in{\cal P}_1$.
These are the points in the sublattice $\Lambda_2$ which are
in the Voronoi set $V_0$ of $\Lambda_s$
when translated to be centered at
$\lambda_1\in{\cal P}_1$.
This ensures that the edge length
$\|\alpha_2(\lambda)-\alpha_1(\lambda)\|^2$ will be
minimized (see Property 3).
We will show that each member of ${\cal L}_1(\lambda_1)$ lies
in a different coset with respect to the sublattice shifts in $\Lambda_s$
(Property 2).
Similarly, we determine the set 
\begin{equation}\label{EqP4}
{\cal L}_2(\lambda_2)=\{\lambda_1\in\Lambda_1:
\lambda_1\in V_0+\lambda_2\}
\end{equation}
for all $\lambda_2\in{\cal P}_2$.
The set of edges emanating from $V_0$ is given by 
\begin{equation}\label{EqP5}
{\cal E}_{edges} = 
\{(\lambda_1,\lambda_2):\lambda_1\in{\cal P}_1~ , \,\,\lambda_2\in
{\cal L}_1(\lambda_1)\} \,\,\bigcup \,\, \{(\lambda_2, \lambda_1): \lambda_2\in{\cal P}_2 ~, \,\,
\,\,\lambda_1\in{\cal L}_2(\lambda_2)\}\,.
\end{equation}
We find a set of coset representatives ${\cal E}_0$ for the equivalence classes of ${\cal E}_{edges}$ modulo $\Lambda_s$.
Property 6 will establish that we can write ${\cal E}_0$ either as
\begin{equation}\label{EqP6}
{\cal E}_0=\{(\lambda_1,\lambda_2):\lambda_1\in{\cal P}_1 \,\,{\rm and} \,\,\lambda_2\in
{\cal L}_1(\lambda_1)\}
\end{equation}
or equally well as
\begin{equation}\label{EqP7}
{\cal E}_0=\{\lambda_2\in{\cal P}_2 \,\,{\rm and} 
\,\,\lambda_1\in{\cal L}_2(\lambda_2)\}\,.
\end{equation}

\item
Matching the edges to the lattice points in the Voronoi set 
is now a straightforward and easily solved assignment problem
(cf. \cite{LuenbergerBlueBook84}).
To formulate this assignment problem we compute the cost given by (\ref{eq:opt_cr})
for each lattice point and
each equivalence class of edges modulo $\Lambda_s$ (taking the minimum over the edge class).
We use only one member from each edge class modulo $\Lambda_s$ in order for the shift invariance property to be satisfied.
This allows us to construct the set of edges which will later be used to label
the points in $V_{\Lambda: \Lambda_s}$.

If there exists a sublattice $\Lambda_{lcm}$ (as defined in Section II) which
is also a geometrically strictly similar, clean sublattice of $\Lambda_1$ and
$\Lambda_2$ the computational complexity of the design can be further
reduced.
For then we need only label the coset representatives for ${\mathcal E}_0$ modulo $\Lambda_{lcm}$.
We will show that this does not reduce the performance of the quantizer 
--- see Property \ref{prop:prop9}.
In this case
we replace the sets ${\cal P}_1$ and ${\cal P}_2$
by the sets
${\cal P}_1^{'}=V_{\Lambda_{1}:\Lambda_{lcm}}(0)$ and
${\cal P}_2^{'}=V_{\Lambda_{2}:\Lambda_{lcm}}(0)$.
The rest of the procedure is unchanged.

\end{enumerate}

\subsection{Properties of the quantizer}
In this section we state some of the properties of the construction
proposed in Section \ref{subsec:lattice}.
We have imposed the following
restrictions on the labeling scheme:

\begin{enumerate}
\item The labels are invariant under shifts by the product sublattice $\Lambda_s$.
\item The labels arise from different cosets of the product sublattice:
i.e.
if $(\lambda_1,\lambda_2)$ and $(\lambda_1, \lambda^{'}_2 )$ are valid edges
then $\lambda_2$ and $\lambda^{'}_2$ are in different
cosets with respect to the product sublattice.
\end{enumerate}

\begin{property}
\label{claim:claim0}
Each member of ${\cal L}_1(\lambda_1)$ lies
in a different coset with respect to the sublattice shifts in $\Lambda_s$,
and $|{\cal L}_1(\lambda_1)|=N_1$.
Similarly each member of ${\cal L}_2(\lambda_2)$ lies
in a different coset with respect to the sublattice shifts in $\Lambda_s$,
and $|{\cal L}_2(\lambda_2)|=N_2$.
\end{property}
\begin{proof}
Let $\lambda_2,\lambda_2^{'}\in{\cal L}_1(\lambda_1)$, and
$\lambda_2^{'}=\lambda_2+\lambda_s$ for some $\lambda_s\in \Lambda_s$.
Then $\lambda_2^{'}-\lambda_1=\lambda_2-\lambda_1+\lambda_s$.
Hence $\lambda_2-\lambda_1$ and $\lambda_2^{'}-\lambda_1$ cannot both 
lie in $V_0$.
But since $\lambda_2,\lambda_2^{'}\in{\cal L}_1(\lambda_1)$,
$\lambda_2-\lambda_1$
and $\lambda_2^{'}-\lambda_1$ are in $V_0$,
a contradiction.
Thus each $\lambda_2\in{\cal L}_1(\lambda_1)$ is in a different coset 
with respect to the sublattice shifts in $\Lambda_s$.
Now $\{V_0+\lambda_1^{'}\}_{\lambda_1^{'}=\lambda_1+\lambda_s,\lambda_s\in
\Lambda_s}$ is a partitioning of the points of $\Lambda$, and each of these disjoint sets contains points
from different cosets of $\Lambda_2$ (with respect to shifts in $\Lambda_s$).
Since there are only $N_1$ different cosets of $\Lambda_2$,
$|{\cal L}_1(\lambda_1)|\leq N_1$. 
In fact equality must hold, because the space is tiled by such sets
and if there were a $\lambda_1^{'}$ for which $|{\cal L}_1(\lambda_1^{'})|< N_1$
then would be a $\lambda_1$ for which $|{\cal L}_1(\lambda_1)|> N_1$,
which impossible.
An identical proof holds for ${\cal L}_2(\lambda_2)$. 
\end{proof}

\begin{property}
\label{claim:claim1}
${\cal L}_1(\lambda_1)$ consists of the $N_1$ points $\la_2 \in \Lambda_2$
closest to $\lambda_1$ subject to the constraint that each $\lambda_2$ is
in a different coset.
\end{property}
\begin{proof}
We know that $\lambda_2\in{\cal L}_1(\lambda_1) \Leftrightarrow \lambda_2\in
V_0+\lambda_1  \Leftrightarrow \lambda_2-\lambda_1\in V_0$.
Then $\lambda_2-\lambda_1\in V_0 \Leftrightarrow ||\lambda_2-\lambda_1||^2
\leq ||\lambda_2-\lambda_1 +\lambda_2||^2$, for all $\lambda_s\in \Lambda_s$.
Thus for any $\lambda_2^{'} = \lambda_2+\lambda_s, \lambda_s\neq 0$
we have $||\lambda_2-\lambda_1||^2\leq ||\lambda_2^{'}-\lambda_1||^2$,
and the claim follows.
\end{proof}

\begin{property}
\label{claim:claim2}
$\lambda_2\in {\cal L}_1(\lambda_1)  \Leftrightarrow \lambda_1 \in {\cal L}_2(\lambda_2)$.
\end{property}
\begin{proof}
For clean lattices, if $x\in V_0$ then $-x\in V_0$ \cite{SPLAG}.
Then $\lambda_2\in {\cal L}_1(\lambda_1)\Leftrightarrow \lambda_2\in V_0+\lambda_1
\Leftrightarrow \lambda_2-\lambda_1 \in V_0 \Leftrightarrow \lambda_1-\lambda_2
\in V_0 \Leftrightarrow \lambda_1 \in V_0+\lambda_2 \Leftrightarrow
\lambda_1 \in {\cal L}_2(\lambda_2)$.
\end{proof}

\begin{property}
\label{prop:prop0}
As lattice points in $\Lambda$ are labeled,
the number of times each point from $\Lambda_1$ is used is $N_1$ and
the number of times each point from $\Lambda_2$ is used is $N_2$.
\end{property}
\begin{proof}
Let $N(\lambda_1)$ denote the number of lattice points 
labeled by
$\lambda_1\in \Lambda_1$.
Certainly $N(\lambda_1)\geq N_1$
since $\lambda_1$ is used $N_1$ times when we form the edges 
$\{(\lambda_1,\lambda_2):\lambda_2\in {\cal L}_1(\lambda_1)\}$.
If $\lambda_1$ is used in more than $N_1$ labels then there is a valid edge
$(\lambda_1,\lambda_2)$ with
$\lambda_2\notin {\cal L}_1(\lambda_1)$.
But this is impossible by Property \ref{claim:claim2}.
Therefore $N(\lambda_1)=N_1$ for all $\lambda_1\in \Lambda_1$ and similarly
$N(\lambda_2)=N_2$ for all $\lambda_2\in \Lambda_2$.
\end{proof}

\begin{property}
\label{prop:prop1}
The number of cosets in
the edge set ${\cal E}_0$
modulo $\Lambda_s$ is equal to the 
number of lattice points in $V_0$.
\end{property}
\begin{proof}
Consider the edge set ${\cal E}_0 = \{(\lambda_1,\lambda_2):
\lambda_1\in {\cal P}_1, \lambda_2\in{\cal L}_1(\lambda_1)\}$.
From Property \ref{claim:claim0} each $\lambda_2\in{\cal L}_1(\lambda_1)$ lies
in a different coset modulo $\Lambda_s$ and hence
each edge $(\lambda_1,\lambda_2)\in{\cal E}_0^{(1)}$ lies in a different
coset. As $|{\cal E}_0^{(1)}| = N_1N_2$, there are at least that
many cosets in the edge set.
\end{proof}

\begin{property}
\label{prop:prop2}
The labeling scheme produces a unique label for each lattice point.
\end{property}
\begin{proof}
This is immediate from the fact that the labels for the cosets of $\Lambda / \Lambda_s$ are taken from distinct cosets of ${\mathcal E}_0 / \Lambda_s$.
\end{proof}

\begin{property}
\label{prop:prop3}
The labeling scheme minimizes the cost criterion
given in {\rm (\ref{eq:Lag})} subject to the coset restriction.
\end{property}
\begin{proof}
This is an immediate consequence of Property \ref{claim:claim1}.
\end{proof}

\begin{property}
\label{prop:prop9}
Suppose $N_1$ and $N_2$ are not relatively prime, and there exists a sublattice $\Lambda_{lcm}$ with index $lcm \{ N_1, N_2 \}$ in $\Lambda$ which is a geometrically strictly similar,
clean sublattice of $\Lambda_1$ and $\Lambda_2$, and contains $\Lambda_s$.
Then we may construct the labeling to be invariant under shifts by
$\Lambda_{lcm}$, and obtain the same edge set as if we used the product lattice $\Lambda_s$.
With this procedure it is necessary to label only $lcm \{ N_1, N_2 \}$ lattice
points rather than $N_1 N_2$ points.
\end{property}
\begin{proof}
If such a $\Lambda_{lcm}$ exists then we just need to show that
the edge set constructed by using the algorithm with $\Lambda_s$ 
can be produced by sublattice shifts of the edge set constructed using
$\Lambda_{lcm}$.
As we saw in the proof of Property \ref{prop:prop1}, the coset 
representatives for
the edge set are constructed by using
${\cal E}_0=\{(\lambda_1,\lambda_2):
\lambda_1\in {\cal P}_1, \lambda_2\in{\cal L}_1(\lambda_1)\}$.
However, ${\cal P}_1 = \bigcup_{\lambda_{lcm}\in\Lambda_{lcm}}
\bigcup_{\lambda_1
\in{\cal P}_1^{'}} (\lambda_1+\lambda_{lcm})$.
Therefore ${\cal E}_0 = \bigcup_{\lambda_{lcm}\in\Lambda_{lcm}}
{\cal E}_0^{'}$, where ${\cal E}_0^{'}=\{(\lambda_1,\lambda_2):
\lambda_1\in {\cal P}_1^{'}, \lambda_2\in{\cal L}_1(\lambda_1)\}$.
It follows that there are exactly $lcm \{ N_1,N_2 \}$ coset leader edges in
${\cal E}_0$ with respect to the sublattice $\Lambda_{lcm}$ and
they are given in ${\cal E}_0^{'}$.
Therefore, by matching the cosets of the edges modulo
$\Lambda_{lcm}$ with the lattice points in the Voronoi set for
$\Lambda_{lcm}$, using the assignment algorithm (as before),
and then shifting by $\Lambda_{lcm}$ we produce exactly
the same labeling as we obtained using $\Lambda_s$.
\end{proof}

\begin{property}
\label{prop:prop7}
If there exist several labeling schemes
achieving the same cost we can mix these configurations to achieve
different levels of asymmetry.
A sufficient condition for this to occur is for the number of
unique representation points to be smaller than the number of lattice 
points in the product lattice.
\end{property}
\begin{proof}
The number of representation points is equal to the number of lattice
points in the Voronoi set $V_0$ (see Property \ref{prop:prop1}).
Therefore, if there are some representation points which overlap ({\em i.e.}
the number of unique representation points is less than the number
of points in $V_0$, then there is more than one labeling scheme that produces
the same Lagrangian $\gamma_1\bar{d}_1+\gamma_2\bar{d}_2$, with each
labeling producing different $\bar{d}_1,\bar{d}_2$.
Suppose one extremal configuration produces the lowest $\bar{d}_1^{min}$
and therefore produces the largest $\bar{d}_2^{max}$ and if the other
extremal configuration produces, the highest $\bar{d}_1^{max}$ and
the lowest $\bar{d}_2^{min}$.
Then by using the first configuration in proportion $\alpha$ and the
second in proportion $\bar{\alpha}=1-\alpha$ one can produce
side distortions $\bar{d}_1=\alpha\bar{d}_1^{min}+\bar{\alpha}\bar{d}_1^{max}$
and $\bar{d}_2=\alpha\bar{d}_2^{max}+\bar{\alpha}\bar{d}_2^{min}$.
Thus by keeping the Lagrangian cost the same, one can obtain different
levels of asymmetry in the distortions $\bar{d}_1,\bar{d}_2$.
\end{proof}

\subsection{Example}
In this section we illustrate the design procedure with
an example in two dimensions using the lattice 
$\ZZ^2$.
We choose $|\Lambda_1|=5$ and $|\Lambda_2|=9$.
Portions of the two sublattices are shown in Figure \ref{fig:fig1}
where the points of $\Lambda_1$ are marked with circles, the points of $\Lambda_2$ with 
crosses, and the
points of $\Lambda_s$ with both
circles and crosses.\footnote{In the enhanced (pdf) version of this document the circles are blue and the crosses are red.}
There are 45 points in the Voronoi set $V_0$ for $\Lambda_s$.
The set ${\cal P}_1$ contains 9 points of $\Lambda_1$ and the set ${\cal P}_2$ contains 5 points of $\Lambda_2$.
The edges ${\cal E}_{edges}$ (see Eq. (\ref{EqP5})) emanating from the
points of $V_0$ are also shown.
These are found using the
sets ${\cal L}_1(\lambda_1)$ and
${\cal L}_2(\lambda_2)$
for $\Lambda_1 \in \Lambda_1$, $\lambda_2 \in \Lambda_2$.
For example, if we take the point $\lambda_1=(2,1) \in {\cal P}_1$, we see that there are 5 points in the set 
${\cal L}_1(\lambda_1)$, namely $\{(0,0),(0,3),(3,3),(6,0),(3,0)\}$.
Note that there are several edges emanating from $V_0$ which
are a sublattice $\Lambda_s$ shift apart. 
For example the edge $\{(-2,-1),(-6,0)\}$ is a sublattice $\Lambda_s$
shift away from the edge $\{(4,2),(0,3)\}$.
To satisfy the shift invariance constraint, we must use only
one of these edges to label a point in $V_0$.
This constraint
is built into the optimization procedure.
The result of the optimization procedure is illustrated in Figure \ref{fig:fig2}.
Here we have shown only the points in $\Lambda_0$.
The points in $\Lambda_1 \cap \Lambda_0$ are marked by circles and those in
$\Lambda_2 \cap \Lambda_0$ by crosses.
Each point carries a pair of labels $(\lambda_1, \lambda_2)$ with
$\Lambda_1 \in \Lambda_1$, $\lambda_2 \in \Lambda_2$.
In this example
we have set $\gamma_1=9$ and $\gamma_2=5$, which
determines the respective distortions $\bar{d}_i$ obtained
by the design.
A comparison of these distortions with that predicted by information
theory is given in Section \ref{sec:results}.

\begin{figure}[htb]
\begin{center}
\includegraphics[scale=0.6]{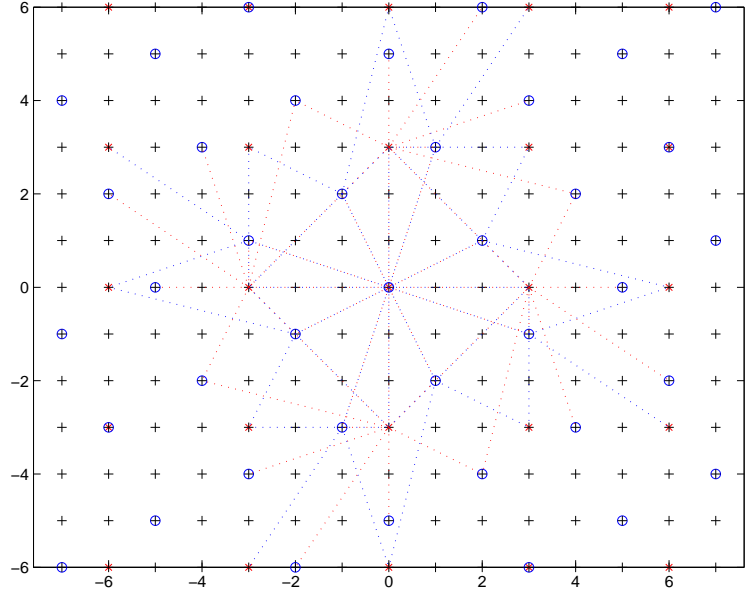}
\end{center}
\caption{Edges emanating from the Voronoi cell of the sublattice.}
\label{fig:fig1}
\end{figure}

\begin{figure}[htb]
\begin{center}
\includegraphics[scale=0.6]{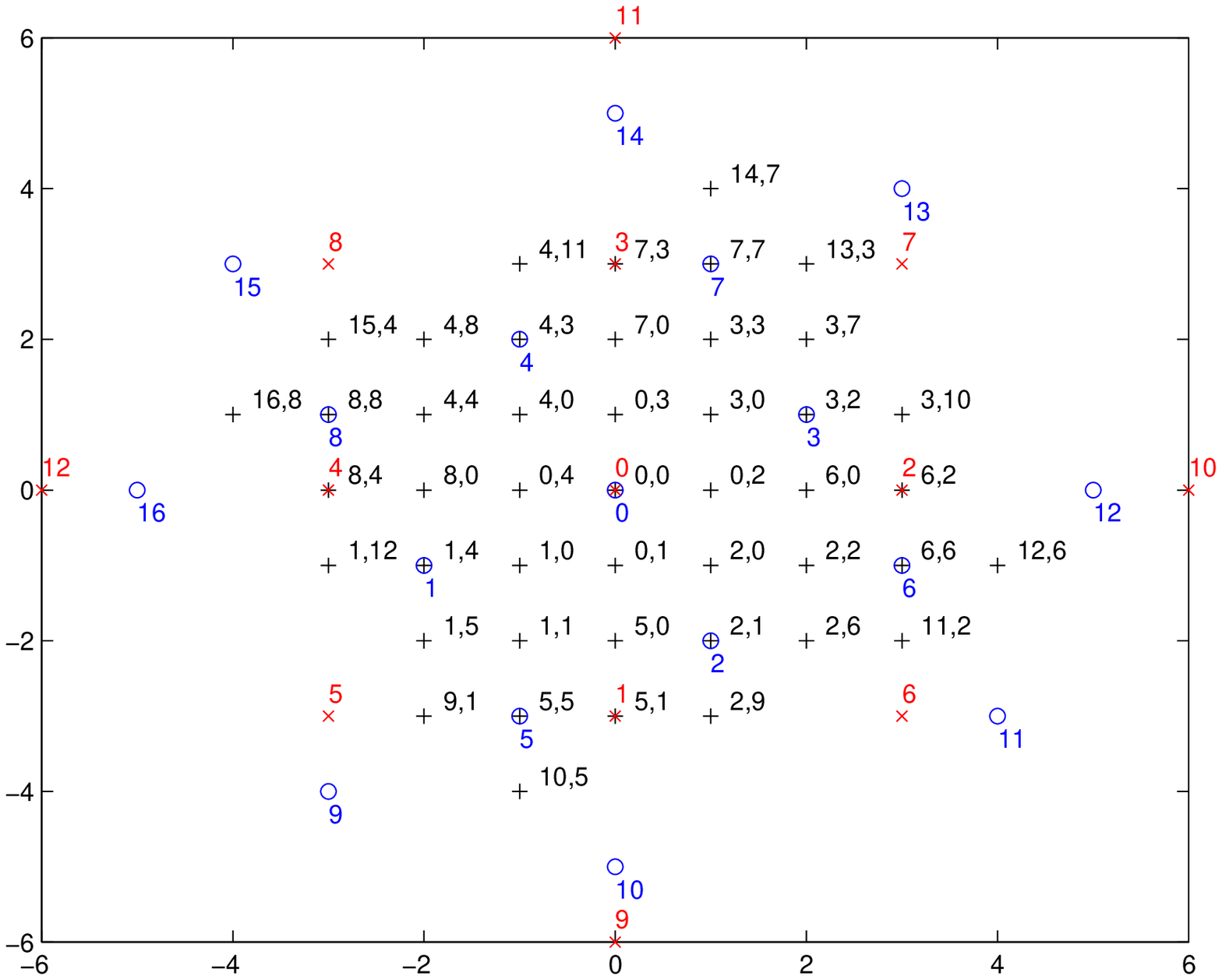}
\end{center}
\caption{Labels generated by the algorithm.}
\label{fig:fig2}
\end{figure}

\section{Good lattices}
\label{sec:good}

The lattices that we will investigate and apply in this paper are
$\ZZ^n$ for $n=1$, 2 or a multiple of 4,
together with the root lattices $D_4$ and $E_8$ \cite{SPLAG}.
The analysis could be extended to treat other lattices such as
$\ZZ^3$, $\ZZ^6$, the 12-dimensional Coxeter-Todd lattice, the 16-dimensional
Barnes-Wall lattice or the 24-dimensional Leech lattice (\cite{SPLAG}, \cite{NS00}), but we shall not discuss these here.

\subsection{The construction of similar sublattices}\label{sec:css}

We begin with the observation that multiplication of points in the square lattice $\ZZ^2$ (regarded as points in the complex plane) by $1+i$ produces a similar
sublattice of index 2.
All our sublattices will be constructed by generalizing this remark.

We will make use of five types of integers:
$\ZZ$, the ordinary {\em rational integers}; ${\mathcal G}$, the ring of
{\em Gaussian integers} $\{a+bi : a,b \in \ZZ \}$, where
$i = \sqrt{-1}$; ${\mathcal J}$, the ring of {\em Eisenstein integers}
$\{a +b \omega : a, b \in \ZZ \}$, where
$\omega = e^{2 \pi i /3}$; ${\mathbb H}_0$, the ring of
{\em Lipschitz integral quaternions}
$\{a+bi+cj+dk: a,b,c,b \in \ZZ \}$, where $i,j, k$ are the familiar unit quaternions;
and ${\mathbb H}_1$, the ring of {\em Hurwitz integral quaternions}
$\{a+ bi + cj + dk : a,b,c,d$ all in $\ZZ$ or all in $\ZZ + \frac{1}{2} \}$.
Other rings of integers could also be used, but these suffice for the lattices
considered in this paper.

If $\Lambda = \ZZ$, multiplication of lattice points by $\xi \in \ZZ$ gives $\xi \ZZ$, a similar sublattice of index $N = |\xi|$.

If $\Lambda = \ZZ^2 = {\mathcal G}$, multiplication by the Gaussian integer
$\xi = a+ bi \in \mathcal G$ gives a similar sublattice $\Lambda ' = \xi \Lambda$ of index
$N = a^2 + b^2$.
A number $N$ is of the form $a^2 + b^2$ if and only if it is of the form
\begin{equation}\label{Eq10f2}
2^{e_1} \prod_{p_i \equiv 1(4)} p_i^{f_i}~
\prod_{q_j \equiv 3(4)} q_j^{2g_j} ~,
\end{equation}
where the first product is over primes $p_i$ congruent to 1 $(\bmod~4)$, the
second product is over primes $q_j$ congruent to 3 $(\bmod~4)$ and
$e_1$, $f_i$ and $g_j$ are nonnegative integers.
These indices are the numbers
\begin{equation}\label{Eq10f1}
1,2,4,5,8,9,10,13,16,17,18,20, \ldots
\end{equation}
(Sequence \htmladdnormallink{A1481}{http://www.research.att.com/cgi-bin/access.cgi/as/njas/sequences/eisA.cgi?Anum=001481}
of \cite{OEIS}).

If $\Lambda = A_2 = {\mathcal J}$, the planar hexagonal lattice, multiplication by
the Eisenstein integer $\xi = a+ b \omega \in \mathcal J$ gives a similar sublattice
$\Lambda ' = \xi \Lambda$ of index $N = a^2 + ab + b^2$.
A number $N$ is of the form $a^2 + ab + b^2$ if and only if
primes congruent to 2 $(\bmod ~3)$ appear to even powers.
These indices are the numbers
\begin{equation}\label{Eq10g1}
1,3,4,7,9,12,13,16,19,21,25,27, \ldots
\end{equation}
(Sequence \htmladdnormallink{A3136}{http://www.research.att.com/cgi-bin/access.cgi/as/njas/sequences/eisA.cgi?Anum=003136} of \cite{OEIS}).

It is shown in \cite{sim} that the above conditions are also necessary:
if $\ZZ, \ZZ^2$ or $A_2$ has a similar sublattice of index $N$ then $N$
must have the form described in the preceding paragraphs.

For the lattices $\Lambda = \ZZ^4$, $\ZZ^8$, $\ZZ^{12} , \ldots, D_4$ and $E_8$ a necessary condition for the existence of a geometrically similar sublattice
of index $N$ is that $N$ should be of the form $m^{L/2}$ for some integer $m$,
where $L$ is the dimension.
This condition is also sufficient, since such sublattices can be obtained
by writing $m= a^2 + b^2 +c^2 + d^2$, regarding $\Lambda$ as a sublattice of
${\mathbb H}_0^{L/4}$, and multiplying $\Lambda$ on the left or on the right by the quaternion $\xi = a + bi + cj + dk$.
Left and right multiplications in general give different sublattices.
In the case of $D_4$ and $E_8$ we may also multiply by Hurwitz integral quaternions to obtain
further similar sublattices.
 
Odd-dimensional lattices of dimension greater than 1 are less
interesting.
For a lattice $\Lambda$ of odd dimension $L$ has a geometrically similar
sublattice of index $N$ if and only if $N$ is an $L$-th power, say $m^L$, and sublattices of this index can be obtained by scalar multiplication of $\Lambda$ by $m$ (see \cite{sim}).

The {\em norm} of a quaternion $\xi = a+ bi + cj +dk$ is $\xi \bar{\xi} = a^2 + b^2 + c^2 + d^2$
where the bar denotes quaternionic conjugation.
If $\xi$ belongs to one of the above rings,
the index of the sublattice $\xi \Lambda$ (or $\Lambda \xi$) in $\Lambda$,
$[\Lambda : \xi \Lambda ]$, is equal to
$(\xi \bar{\xi} )^{L/2}$, where $L$ is the dimension and the bar 
is complex or quaternionic conjugation as appropriate.

\subsection{Clean sublattices}\label{sec:clean}

In dimension one, the sublattice $\xi \ZZ$ is clean if and only if $\xi$ is odd.

Reference \cite{sim} gives necessary and sufficient conditions for a similar sublattice of any two-dimensional lattice to be clean.
In particular, the sublattice $\xi \ZZ^2$ $(\xi = a+ ib )$ is clean if and only if $N = a^2 + b^2$ is odd.
These indices are obtained by setting $e_1 =0$ in \eqn{Eq10f2}:
$$1,5,9, 13, 17,25,29,37,41, 45 , \ldots$$
(Sequence \htmladdnormallink{A57653}{http://www.research.att.com/cgi-bin/access.cgi/as/njas/sequences/eisA.cgi?Anum=057653}).

The sublattice $\xi A_2$ $(\xi = a+b \omega )$ is clean if and only if $a$ and $b$ are relatively prime.
It follows that $A_2$ has a clean similar sublattice of index $N$ if and only if $N$ is a product of primes congruent to 1 $(\bmod~6)$.
These are the numbers
\begin{equation}\label{Eq10n1}
1,7,13,19,31,37,43,49,61,67, \ldots
\end{equation}
(Sequence \htmladdnormallink{A57654}{http://www.research.att.com/cgi-bin/access.cgi/as/njas/sequences/eisA.cgi?Anum=057654}).

The existence of clean sublattices in dimensions greater than 2 was not considered in \cite{sim}.

We can give a fairly complete answer for the lattices $\ZZ^L$, $L \ge 1$.

\begin{theorem}\label{CLZ1}
Suppose $L \ge 1$ and $\ZZ^L$ has a geometrically similar sublattice $\Lambda '$ of index $N$.
Then $\Lambda '$ is clean if and only if $N$ is odd.
\end{theorem}

\paragraph{Proof.}
(If)
Let $\Lambda' = \phi (\ZZ^L )$, where $\phi$ is a similarity, and let $\Lambda '' = \phi^{-1} (\ZZ^L )$.
If $\phi$ multiplies lengths by $c_1$ (as in \eqn{EqU1}) then $N= c_1^L$.
Suppose $N$ (and hence $c_1$) is odd.
Let $\Lambda '$ have generator matrix $K$, with $KK^{tr} = c_1^2 I_L = mI_L$,
where $m= c_1^2$.
(In the notation of \eqn{EqU1}, $U_1 = c_1 K_1 = K$.)
Since $\Lambda '$ is a sublattice of $\ZZ^L$, the entries of $K$ are integers.
Then $\Lambda ''$ has generator matrix $K^{-1} = \frac{1}{m} K^{tr}$.

We must show that there are no points of $\ZZ^L$ on the boundary of the Voronoi cell of $\Lambda '$, or equivalently that there are no points of $\Lambda ''$ on the boundary of the Voronoi cell of $\ZZ^L$.

It is enough to consider just one face of the Voronoi cell of $\ZZ^L$,
say that consisting of the points
$P = \left ( \frac{1}{2} , \frac{x_2}{2} , \frac{x_3}{2} , \ldots, \frac{x_{L-1}}{2} \right)$, where
$|x_i | \le 1$ for $2 \le i \le L-1$.
If $P \in \Lambda ''$ then there is a vector $u = (u_1, \ldots, u_L ) \in \ZZ^L$ such that
\begin{equation}\label{Eq52a}
P = \frac{1}{m} u K^{tr} ~.
\end{equation}
Equating the first components we get that
$$
\frac{1}{2} = \frac{1}{m} ~\mbox{times a vector with integer entries} ~.
$$
Since $m$ is odd this is impossible.

(Only if)
Suppose $N$ (and hence $c_1$ and $m$)
is even.
We claim that all the vertices of the Voronoi cell for $\ZZ^L$ (i.e. all the deep holes in $\ZZ^4$ in the notation of \cite{SPLAG}) belong to $\Lambda ''$.
In fact, \eqn{Eq52a} implies that $u= PK$.
Let $P$ be a vector of the form $(\pm \frac{1}{2}, \pm \frac{1}{2}, \ldots, \pm \frac{1}{2} )$, and let $K = (k_{ij} )$.
From $KK^{tr} = K^{tr} K = mI_L$ we have $\sum_{i=1}^L k_{ij}^2 = m$ and, since
$k_{ij}^2 \equiv k_{ij}$ $(\bmod~2)$, $\sum_{i=1}^L k_{ij}$ is even (for all $j$).
Hence $PK$ has integer entries and is in $\ZZ^L$.
\hfill$\blacksquare$

The following corollary summarizes our results about $\ZZ^L$
for the values of $L$ that we are interested in.
Note that since $\ZZ^L$ has no ``handedness'', there is essentially
no difference between ``similar'' and ``strictly similar'' for this lattice.

\begin{corollary}\label{CLZ2}
$\ZZ^L$ has a geometrically similar sublattice of index $N$ if and only if
\begin{itemize}
\item
$N$ is an $L^{th}$ power, if $L$ is odd
\item
$N$ is of the form $a^2 + b^2$, if $L =2$
\item
$N$ is any integer, if $L = 4k$, $k \ge 1$
\end{itemize}
In each case the sublattice is clean if and only if $N$ is odd.
The same results hold if ``similar'' is replaced by ``strictly similar''.
\end{corollary}

For $D_4$ we have only a partial answer.

\begin{theorem}\label{CLD1}
If $M$ is $7$ or a product of primes congruent to $1$ $(\bmod~4)$ then $D_4$
has a geometrically strictly similar, clean sublattice of index $M^2$.
The values of $M$ mentioned are
\begin{equation}\label{Eq45aa}
1,5,7,13,17,25,29,37,41,53, \ldots
\end{equation}
($7$ together with Sequence
\htmladdnormallink{A4613}{http://www.research.att.com/cgi-bin/access.cgi/as/njas/sequences/eisA.cgi?Anum=004613}).
\end{theorem}

\paragraph{Proof.}
We take our standard version of the $D_4$ lattice to have minimal norm 2
(as in \cite{SPLAG}) and generator matrix
\begin{equation}\label{Eq45a}
G = \left[
\begin{array}{rrrr}
1 & 1 & 0 & 0 \\
-1 & 0 & 1 & 0 \\
0 & -1 & 0 & 1 \\
0 & -1 & 0 & -1
\end{array}
\right] ~.
\end{equation}
The four rows $v_1, v_2, v_3, v_4$ of $G$ correspond to the nodes of the
Coxeter diagram for $D_4$ shown in Fig.~4,
where $v_i \cdot v_i =2$ $(i=1, \ldots, 4)$, two nodes that are joined by an edge correspond to vectors with inner product $-1$, and two nodes that are not joined by an edge are orthogonal.

\begin{figure}[htb]
\begin{center}
\input f45.pstex_t
\end{center}
\vspace*{+.1in}
Figure 4: Coxeter diagram for any lattice that is geometrically similar to $D_4$: there are four generating vectors $v_1,v_2,v_3,v_4$ satisfying $v_1 \cdot v_1 = v_2 \cdot v_2 = v_3 \cdot v_3 = v_4 \cdot v_4$, $v_i \cdot v_j = - \frac{1}{2} v_1 \cdot v_1$ if nodes $v_i$ and $v_j$ are joined by an edge, and $v_i \cdot v_j =0 ~(i \neq j)$ otherwise.
\label{F45}
\end{figure}

We regard $D_4$ as a subset of ${\mathcal H} = \{w + xi + yj + zk : w,x,y,z \in {\mathbb R} \}$, the space of real quaternions.
Our sublattices $\Lambda '$ will be constructed by multiplying $D_4$ either on the left or on the right by appropriate Hurwitzian integral quaternions.
If $\xi = a+ bi + cj + dk \in {\mathcal H}$ then $\xi D_4$ has generator matrix
 $GL_\xi$, where
\begin{equation}\label{Eq51b}
L_\xi = \left[
\begin{array}{rrrr}
a & b & c & d \\
-b & a & d & -c \\
-c & -d & a & b \\
-d & c & -b & a
\end{array}
\right] ~,
\end{equation}
and $D_4 \xi$ has generator matrix $GR_{\xi}$, where
\begin{equation}\label{Eq51c}
R_\xi = \left[
\begin{array}{rrrr}
a & b & c & d \\
-b & a & -d & c \\
- c& d & a & -b \\
-d & -c & b & a
\end{array}
\right] ~.
\end{equation}
Note that
\begin{equation}\label{Eq51d}
L_\xi L_\xi^{tr} = R_\xi R_\xi^{tr} = mI_4 ,~~
L_\xi R_\xi = R_\xi L_\xi ~,
\end{equation}
where $m = \xi \bar{\xi} = a^2 + b^2 + c^2 + d^2$.

We will show that under certain conditions $\xi D_4$ and $D_4 \xi$ are clean sublattices.
We only give the proof for $D_4 \xi$, the other case being completely analogous.

The Voronoi cell for $D_4$ is a 24-cell, with 24 octahedral 
faces \cite{cell}, \cite{SPLAG}.
A typical face (they are all equivalent) is that lying in the hyperplane
\begin{equation}\label{Eq51e}
X \cdot v_1 = \frac{1}{2} v_1 \cdot v_1 ~,
\end{equation}
having center $\delta_0 = \frac{1}{2} v_1$ and six vertices
\begin{eqnarray}\label{Eq51f}
\delta_1 & = & \frac{1}{2} (2v_1 + v_3 + v_4 ) \,, \nonumber \\
\delta_2 & = & \frac{1}{2} (-v_3 - v_4) \,, \nonumber \\
\delta_3 & = & \frac{1}{2} (2v_1 + v_2 + v_3 ) \,, \nonumber \\
\delta_4 & = & \frac{1}{2} (-v_2 - v_3) \,, \nonumber \\
\delta_5 & = & \frac{1}{2} (2v_1 + v_2 + v_4 ) \,, \nonumber \\
\delta_6 & = & \frac{1}{2} (-v_2 - v_4 )
\end{eqnarray}
(see Fig. \ref{F46}).
\setcounter{figure}{4}
\begin{figure}[htb]
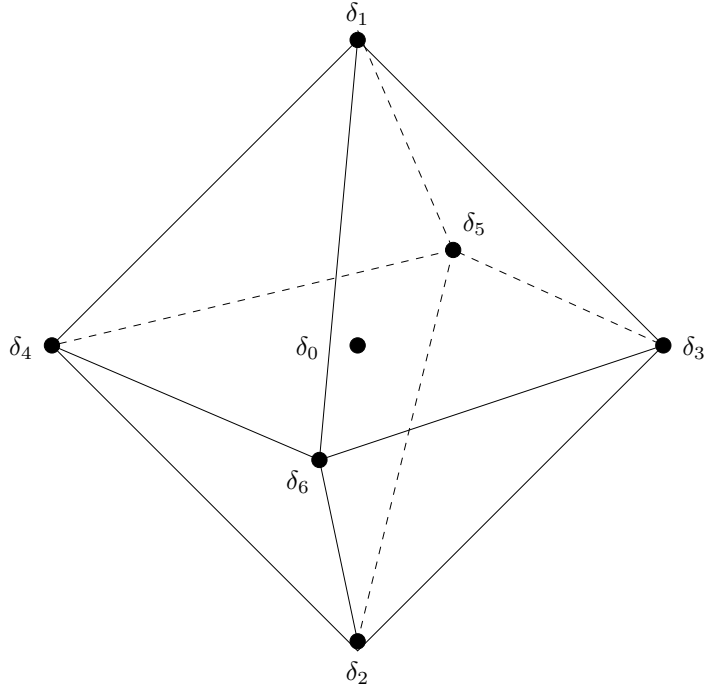

\begin{center}
\input f46.pstex_t
\end{center}
\caption{Labeling for center and vertices of octahedral face of Voronoi cell for $D_4$.}
\label{F46}
\end{figure}
A point
$X$ belongs to this face if and only if it satisfies (\ref{Eq51e}) and
\begin{equation}\label{Eq51g}
|(X- \delta_0 ) \cdot ( \delta_1 - \delta_0 ) | +
| (X- \delta_0 ) \cdot (\delta_3 - \delta_0 ) | +
| (X- \delta_0 ) \cdot (\delta_5 - \delta_0)| \le
\frac{1}{4} v_1 \cdot v_1 \,.
\end{equation}

Let $\Lambda' = D_4 \xi$, where $\xi$ is a quaternion of the form
\begin{equation}\label{Eq52ap}
\xi = \frac{\alpha}{2} + \frac{\alpha}{2} i + \frac{\beta}{2} j + \frac{\beta}{2} k ~,
\end{equation}
and $\alpha$ and $\beta$ are odd, positive, relatively prime integers.
The norm of $\xi$ is $\frac{1}{2} (\alpha^2 + \beta^2 )$.
Then we claim that $\Lambda '$ is clean.

To show this, we begin by computing the generator matrix for $\Lambda '$:
\begin{eqnarray}\label{Eq52b}
G' & = & GR_\xi \nonumber \\
& = & \left[ \begin{array}{cccc}
0 & \alpha & 0 & \beta \\ [+.1in]
\frac{-\alpha - \beta}{2} & \frac{- \alpha + \beta}{2} &
\frac{\alpha - \beta}{2} & \frac{- \alpha - \beta}{2} \\ [+.2in]
\frac{\alpha - \beta}{2} & \frac{- \alpha - \beta}{2} &
\frac{\alpha+ \beta}{2} & \frac{\alpha - \beta}{2} \\ [+.2in]
\frac{\alpha + \beta}{2} & \frac{- \alpha + \beta}{2} &
\frac{- \alpha + \beta}{2} & \frac{-\alpha - \beta}{2}
\end{array}
\right] \,,
\end{eqnarray}
and denote its rows by $v'_1$, $v'_2$, $v'_3$, $v'_4$.
We will similarly use primes to denote the center $(\delta'_0 )$ and vertices
$(\delta'_1, \ldots, \delta'_6 )$ of an octahedral face of the Voronoi cell of $\Lambda '$.
From (\ref{Eq51f}) we find that
\begin{eqnarray*}
\delta'_0 & = & \frac{1}{2} (0, \alpha, 0, \beta ) \,, \\
\delta'_1 & = & \frac{1}{2} (\alpha , \alpha , \beta , \beta ) \,, \\
\delta'_3 & = & \frac{1}{2} ( - \beta , \alpha, \alpha , \beta ) \,, \\
\delta'_5 & = & \frac{1}{2} ( 0, \alpha + \beta , 0, - \alpha + \beta ) \,.
\end{eqnarray*}
We must show that it is impossible for a point $X = (w,x,y,z) \in D_4$
to satisfy the primed versions of (\ref{Eq51e}) and (\ref{Eq51g}), which are
\begin{equation}\label{Eq52c}
\alpha x + \beta z = \frac{1}{2} ( \alpha^2 + \beta^2 ) \,,
\end{equation}
\begin{equation}\label{Eq52d}
|\alpha w + \beta y | + |
- \beta w + \alpha y | + | \beta x - \alpha z | \le
\frac{1}{2} (\alpha^2 + \beta^2 ) \,.
\end{equation}
Suppose on the contrary that $(w,x,y,z) \in D_4$ satisfies (\ref{Eq52c}) and
(\ref{Eq52d}).
From (\ref{Eq52c}) we have
\begin{equation}\label{Eq52e}
z = \frac{1}{2 \beta} (\alpha^2 + \beta^2 - \alpha x )
\end{equation}
and from (\ref{Eq52d})
$$| \beta x - \alpha z | \le \frac{1}{2} ( \alpha^2 + \beta^2 ) \,,
$$
which together imply
$$\frac{1}{2} ( \alpha - \beta ) \le x \le \frac{1}{2}
( \alpha + \beta ) .
$$
So we may write $x = \frac{1}{2} ( \alpha + \mu )$, say,
where $\mu$ is an odd integer satisfying $- \beta \le \mu \le \beta$,
and from (\ref{Eq52e})
$$z = \frac{\beta^2 - \alpha \mu}{2 \beta} ~,$$
which implies $\alpha \mu \equiv \beta^2$ $(\bmod~2\beta )$.
Since $\beta$ is odd, $\beta^2 \equiv \beta$ $(\bmod~2\beta )$, and we conclude that
\begin{equation}\label{Eq53a}
\alpha \mu \equiv \beta ~ ( \bmod~2\beta ) \,.
\end{equation}
Thus for some integer $k$, $\alpha \mu - \beta = 2k \beta$, and since $\alpha$ and $\beta$
are relatively prime, $\beta$ must divide $\mu$.
Therefore $\mu = \pm \beta$.
But this is impossible.
For if $\mu = \beta$, $x= \frac{1}{2} (\alpha + \beta )$, $z = \frac{1}{2} (- \alpha + \beta )$,
$\beta x - \alpha z = \frac{1}{2} (\alpha^2 + \beta^2 )$, and then (\ref{Eq52d})
implies $w = y =0$, so $w + x+y+z = \beta \not\in D_4$,
since $\beta$ is odd. A similar
argument applies if $\mu = - \beta$.

So far we have shown that if $\alpha$ and $\beta$ are odd, positive and relatively prime, then the sublattice $D_4 \xi$ is clean,
where $\xi$ is given by (\ref{Eq52ap}).
Suppose $M$ is a product of primes congruent to 1 $(\bmod~4)$.
From the classical theory of quadratic forms (see for example \cite{Cox}), we know that $M= p^2 + q^2$ with $p$ even, $q$ odd and $gcd (p,q) =1$.
We now simply set $\alpha = p+q$ and $\beta = | p-q |$.

It remains to discuss the case $M=7$.
For this we can multiply on the left or on the right by either of the quaternions
$$\xi = \frac{1}{2} + \frac{1}{2} i + \frac{1}{2} j + \frac{5}{2} k ~~\mbox{or}~~
\frac{1}{2} + \frac{3}{2} i + \frac{3}{2} j + \frac{3}{2} k ~.
$$
We omit the straightforward verification that these sublattices are clean.\hfill$\blacksquare$

In the other direction we have:
\begin{theorem}\label{CLD2}
$D_4$ has no clean, geometrically similar sublattice of index $M^2$ if $M$ is $3$, $9$ or $11$.
\end{theorem}

\paragraph{Proof.}
The proof is by exhaustive search, using a computer.
We produced a list of all vectors of norm $2M$ in $D_4$,
and from this we found all similar sublattices of index $M^2$ by finding all sets of four vectors corresponding to the Coxeter diagram of Fig.~4.
Given a sublattice $\Lambda'$, we compute the equations defining an octahedral face of the Voronoi cell from (\ref{Eq51e}) and (\ref{Eq51g}).
Then AMPL \cite{ampl} and CPLEX \cite{cplex} were used to verify that in every case there was a point of $D_4$
on the face.\hfill$\blacksquare$

The preceding discussion has shown that the lattices
$\ZZ$, $\ZZ^2$, $\ZZ^{4k}$ for $k \ge 1$ and $D_4$ have a plentiful supply of clean, geometrically similar sublattices.
We expect the same will be true of the $E_8$ lattice, but this question is presently under investigation.

Finally, we remark that if  $\Lambda'$ is a clean sublattice of $\Lambda$
and $\Lambda''$ is a clean sublattice of $\Lambda'$,
then $\Lambda''$ is a clean sublattice of $\Lambda$.

\subsection{Common sublattices of $\Lambda_1$ and $\Lambda_2$}\label{sec:common}
We begin with a general
comment.
Let $\Lambda_1$ and $\Lambda_2$ be any two sublattices of a lattice $\Lambda$
(they must have the same dimension as $\Lambda$ but are otherwise arbitrary).
Then we may form their intersection
$\Lambda_\cap = \Lambda_1 \cap \Lambda_2$ and their join
$\Lambda_\cup = \Lambda_1 \cup \Lambda_2$, as shown in Fig.~\ref{F47}.
The join is the lattice generated by the vectors of both $\Lambda_1$ and $\Lambda_2$
(and in general is not simply their union).
\begin{figure}[htb]
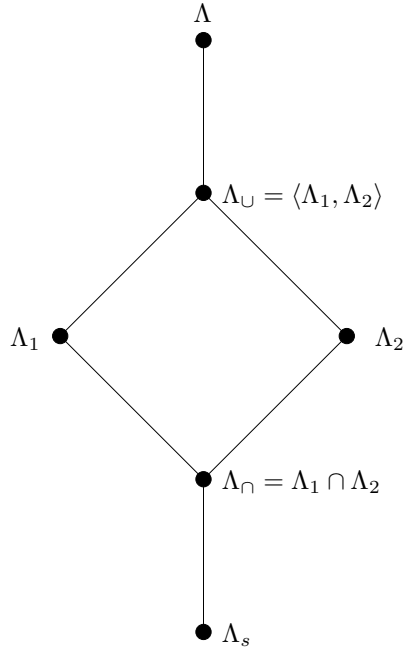

\begin{center}
\input f47.pstex_t
\end{center}
\caption{Intersection, join and ``product'' sublattice of two arbitrary sublattices.}
\label{F47}
\end{figure}
From the second isomorphism theorem of group theory (e.g. \cite{Rotman})
the indices and determinants of these lattices are related by
\begin{equation}\label{Eq54a}
[\La_\cup : \La_1 ] =
[ \La_2 : \La_\cap ] ,~
[\La_\cup : \La_2 ] = [ \La_1 : \La_\cap ] \,,
\end{equation}
\begin{equation}\label{Eq54b}
\det \La_1 ~~\det \La_2 = \det \La_\cup ~~\det \La_\cap ~.
\end{equation}
There are now in general many ways to find a
``product'' sublattice
$\Lambda_s \subset \Lambda_\cap$ with
\begin{equation}\label{Eq34a}
[\Lambda: \Lambda_s] = [ \Lambda : \Lambda_1 ] [ \Lambda : \Lambda_2 ] \,.
\end{equation}

Let $\La$ be one of ${\mathbb Z}$, ${\mathbb Z}^2$ or $A_2$, and let $\La_1 = \xi _1 \La$, $\La_2 = \xi_2 \La$ be geometrically strictly similar sublattices obtained by
multiplying $\La$ by elements of ${\mathbb Z}$, ${\mathcal G}$ or ${\mathcal J}$
respectively.
Since these three rings are unique factorization rings, the notions of greatest common divisor (gcd) and least common multiple (lcm) are
well-defined.
We set $\xi_\cup = gcd (\xi_1, \xi_2 )$, $\xi_\cap = lcm \{ \xi_1, \xi_2 \}$,
and then it is easy to see that $\La_\cup = \xi_\cup \, \La$,
$\La_\cap = \xi_\cap \, \La$.
We can also form the product sublattice 
$\La_s = \xi_1 \xi_2 \, \La$ (see Fig.~\ref{F48}).
The indices of these lattices are given by
\begin{eqnarray}\label{Eq54c}
~[\La : \La_1 ] & = & (\xi_1 \bar{\xi}_1 )^{L/2} , 
~[\La : \La_2 ] = (\xi_2 \bar{\xi}_2)^{L/2} \,, \nonumber \\
~[\La : \La_\cup ] & = & (\xi_\cup \bar{\xi}_\cup )^{L/2} , 
~[\La : \La_\cap ] = (\xi_\cap \bar{\xi}_\cap )^{L/2} \,, \nonumber \\
~[\La : \La_s ] & = & [\La : \La_1 ] [ \La: \La_2 ] \,.
\end{eqnarray}

In dimension $L=1$, (\ref{Eq54c}) implies that

\begin{equation}\label{Eq54cc}
[\La : \La_\cap ] ~ = ~ lcm \{ ~ [\La : \La_1 ], ~ [\La : \La_2 ] ~ \} ~,
\end{equation}
and we can take $\La_{lcm} = \La_\cap$.
However, if $L = 2$, (\ref{Eq54cc}) does not hold in general.

In dimensions 1 or 2, if $\xi_1$ and $\xi_2$ are relatively prime (meaning
$gcd (\xi_1, \xi_2 ) =1$),
we have
$\xi_\cup =1$, $\xi_\cap = \xi_1 \xi_2$, $\La= \La_\cup$,
$\La_s = \La_\cap$.
\begin{figure}[htb]
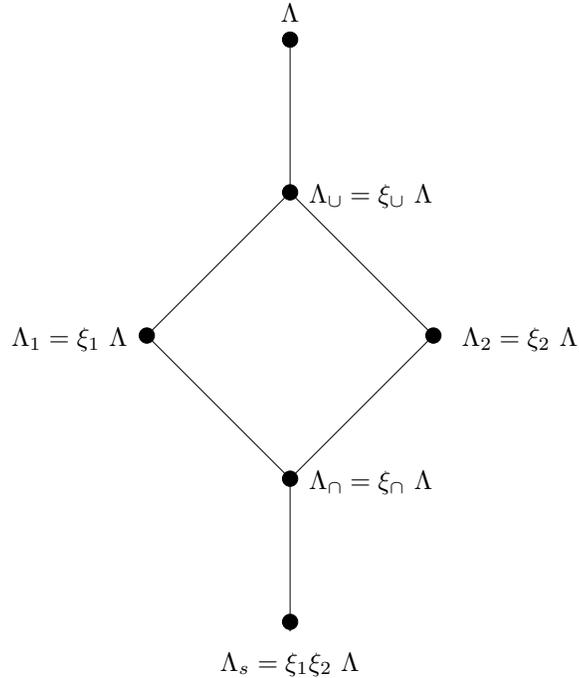

\begin{center}
\input f48.pstex_t
\end{center}
\caption{Join $\La_\cup$, intersection $\La_\cap$ and product $\La_s$ of two sublattices $\La_1$, $\La_2$ of $\La$, where $\La$ is one of ${\mathbb Z}$, ${\mathbb Z}^2$ or $A_2$.}
\label{F48}
\end{figure}

Because the quaternions form a noncommutative ring their arithmetic theory is more complicated.
For example, it is necessary to distinguish between left gcd's and right gcd's.
Both are well-defined in ${\mathbb H}_1$ and also in ${\mathbb H}_0$ as long as at least one of the quaternions involved has odd norm \cite{Dickson}, \cite{HW}.
We plan to discuss this theory and its applications to the study of
sublattices of ${\mathbb Z}^4$ and $D_4$ elsewhere.
In the present paper we will restrict our attention to a narrow class of sublattices, which however will be general enough to provide an adequate
supply of sublattices for our applications.

For ${\mathbb Z}^4$
we choose two Lipschitz integral quaternions $\xi_1$, $\xi_2 \in {\mathbb H}_0$ whose norms are odd and relatively prime.
For $D_4$ we choose two Hurwitz integral quaternions
\begin{eqnarray}\label{Eq54d}
\xi_1 & = & \frac{1}{2} \alpha_1 (1+i) + \frac{1}{2} \beta_1 (j+k) \in {\mathbb H}_1 \,, \nonumber \\
\xi_2 & = & \frac{1}{2} \alpha_2 (1+i) + \frac{1}{2} \beta_2 (j+k) \in {\mathbb H}_1 \,,
\end{eqnarray}
where $\alpha_1$, $\alpha_2$, $\beta_1$, $\beta_2$ are odd positive integers
with $gcd (\alpha_1, \beta_1 ) = gcd (\alpha_2, \beta_2 ) = gcd ((\alpha_1^2 + \beta_2 )/2$, $(\alpha_2^2 + \beta_2^2)/2) =1$.

In both cases we take $\La_1 = \xi_1 \La$, $\La_2 = \La \xi_2$
and $\La_s = \La_\cap = \xi_1 \La \xi_2$ (see Fig.~\ref{F50}).
Then
\begin{eqnarray}\label{Eq54e}
~[\La : \La_1 ] & = & (\xi_1 \bar{\xi}_1 )^2, 
~[\La: \La_2 ] = (\xi_2 \bar{\xi}_2 )^2 \,, \nonumber \\
~[\La : \La_s ] & = & [\La : \La_1 ] [ \La : \La_2 ] \,.
\end{eqnarray}
\begin{figure}[htb]
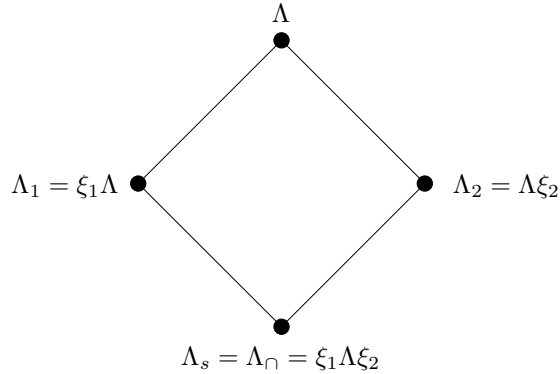

\begin{center}
\input f50.pstex_t
\end{center}
\caption{$\La_1$ (resp. $\La_2$) obtained by multiplying $\La = {\mathbb Z}^4$ or $D_4$ on the left (resp. right) by a quaternion $\xi_1$ (resp. $\xi_2$).}
\label{F50}
\end{figure}

For ${\mathbb Z}^4$ this gives sublattices $\La_1$, $\La_2$ of indices $M_1^2$, $M_2^2$, where $M_1$ and $M_2$ are any two relatively prime odd numbers
(from Corollary \ref{CLZ2}).
For $D_4$, $M_1$ and $M_2$
are any two relatively prime numbers from (\ref{Eq45aa}).

\section{High rate asymptotics: Details}
\label{sec:highrate}

In this section we analyze the distortion of the asymmetric 
multiple description lattice quantizer at high rates.


Let $\Lambda$ be an $L$-dimensional lattice with geometrically
strictly similar, clean
sublattices $\La_1$, $\La_2$, $\La_\cap = \La_1 \cap \La_2$, $\La_s$
(as in Figure ~\ref{F47}),
with indices $N_1,N_2,N_{\cap}$ and $N_s$,
respectively, where $N_s=N_{1}N_2$. It is assumed that $\nu_{\Lambda}$, the volume of a fundamental region for
$\Lambda$, is equal to unity.  A sequence of lattices is then  obtained from the base set of lattices by  scaling each component.  Let $\Lambda_1(n)=n\Lambda_1$, $\Lambda_2(n)=n \Lambda_2$,
$\Lambda_{\cap}(n)=n\Lambda_{\cap}$ and $\Lambda_s(n)=n^{2}\Lambda_s$.
These have indices $N_{\cap}(n)=n^LN_{\cap}$, $N_s(n)=n^{2L} N_{s}$ and  $N_s(n)=N_1(n)N_2(n)$.

We analyze the  rate-distortion performance for the set of lattices
$\{\Lambda,\Lambda_1(n),\Lambda_2(n)$,$\Lambda_{\cap}(n),\Lambda_s(n)\}$.
However, in order to keep the notation simple, we will only use the 
sequence index $n$ when it is necessary to avoid confusion. 
Thus we will write $\Lambda_s$ instead of $\Lambda_s(n)$, 
$N_s$ instead of $N_s(n)$ and so on.
\marginpar{\tiny Need to redefine the Voronoi sets in a consistent way usinng notation that is
manageable.}

Referring to (\ref{eq:Lag}), let 
\begin{equation}
\label{eq:Js}
J_{s}=\sum_{\lambda \in \Lambda} P(\lambda)
\sum_i \gamma_i\|\lambda-\alpha_i(\lambda) \|^2.
\end{equation}
We investigate the high-rate behavior of $J_{s}$
and then find the approximation for $\bar{d}_i,\,i=1,2$.
The latter would also allow us to predict the asymmetry in the
distortion behavior of the quantizer.
The reader is referred to Figure \ref{fig:fig4} for the analysis.
\begin{figure}[htb]
\begin{center}
\includegraphics[scale=0.6]{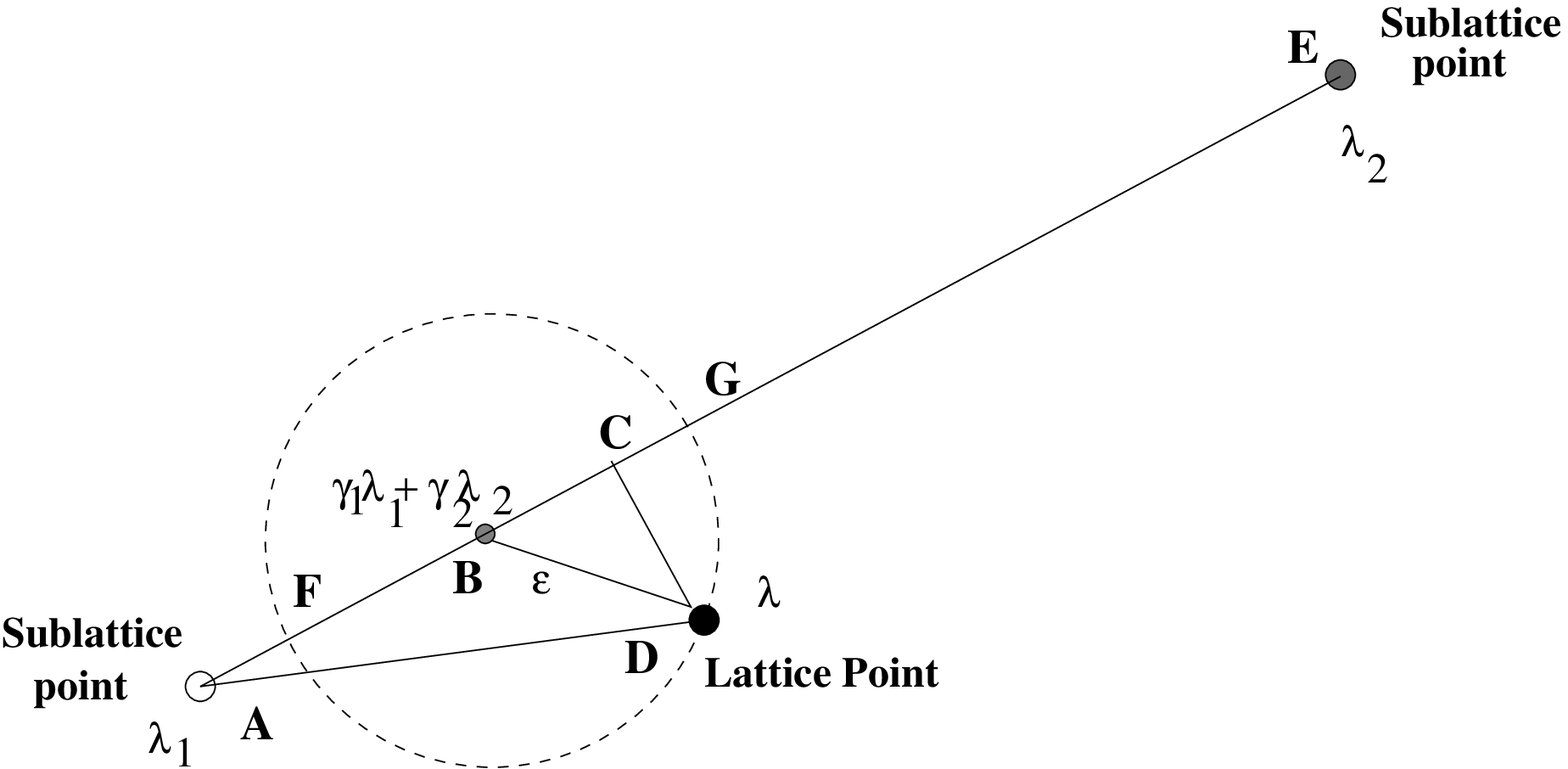}
\end{center}
\caption{Relationship of edgelength and distortion.}
\label{fig:fig4}
\end{figure}
Let 
\begin{equation}
J_{s_1} =\sum_{\lambda \in \Lambda} \frac{\gamma_1\gamma_2}{\gamma_1+\gamma_2}
\|\alpha_2(\lambda)-\alpha_1(\lambda)\|^2 P(\lambda)
\end{equation}
and
\begin{equation}
J_{s_2}=\sum_{\lambda \in \Lambda} (\gamma_1+\gamma_2)\|\lambda-\frac{\gamma_1\alpha_1(\lambda)+\gamma_2
\alpha_2(\lambda)}{\gamma_1+\gamma_2}\|^2 P(\lambda).
\end{equation}
Thus,
\begin{equation}
\label{eqn:JSplit}
J_{s} = J_{s_1} + J_{s_2}.
\end{equation}
Under  the assumption that $\Lambda_{\cap}$ is fine enough for ${\rm P}(\lambda)$ to be considered 
a constant over $V_{\Lambda_{\cap}}$ we obtain
\begin{equation}
\label{eq:lub_Lag}
J_{s_1}=\frac{\gamma_1\gamma_2}{\gamma_1+\gamma_2}\sum_{\lambda'\in\Lambda_{\cap}}{\rm P}(\lambda')
\sum_{\lambda\in V_{\Lambda_\cap}(\lambda')}
||\alpha_1(\lambda)-\alpha_2(\lambda)||^2. 
\end{equation}
By construction, the inner sum in (\ref{eq:lub_Lag}) does not
depend on $\lambda'$. 
Therefore, taking this out of the outer summation and using 
$\sum_{\lambda'\in\Lambda_{\cap}}{\rm P}(\lambda')=1/N_{\cap}$,
we obtain
\begin{equation}
\label{eq:lb_Lag}
{J}_{s_1} = \frac{\gamma_1\gamma_2}{\gamma_1+\gamma_2}
\frac{1}{N_{\cap}} \sum_{\lambda\in V_{\Lambda_\cap}(0)} 
||\alpha_1(\lambda)-\alpha_2(\lambda)||^2,
\end{equation}
which can be written in terms of the edge endpoints
as
\begin{equation}
\label{eq:lb_Lag1}
J_{s_1}=\frac{\gamma_1\gamma_2}{\gamma_1+\gamma_2}\frac{1}{N_{\cap}} \sum_{\lambda_1 \in
V_{\Lambda_{\cap}}(0)} \sum_{\lambda_2 \in V_{\Lambda_s}(\lambda_1)} \|\lambda_1-\lambda_2 \|^2.  
\end{equation}
Here we have used the fact that each edge labels a lattice point and therefore
replaced the sum over the lattice points by the sum over the edges.
Using the Riemann approximation
\begin{equation}
\label{eq:Riemann}
\int_{V_{\Lambda_s(\lambda_1)}}||x - \lambda_1||^2 dx \approx \sum_{\lambda_2\in V_{\Lambda_s(\lambda_1)}} 
||\lambda_2-\lambda_1||^2 \nu_{2}
\end{equation}
for the summation in (\ref{eq:lb_Lag1}) we obtain
\begin{equation}
\label{eq:int_Lag}
{J}_{s_1} = \frac{\gamma_1\gamma_2}{\gamma_1+\gamma_2}
\frac{1}{N_{\cap}}\sum_{\lambda_1\in V_{\Lambda_{\cap}}(0)} \frac{1}{\nu_{2}}
\left [\frac{\int_{V_s(\lambda_1)}||x - \lambda_1||^2 dx}{\nu_s^{(1+2/L)}} \right ]
\nu_s^{(1+2/L)}
\end{equation}
The term within the brackets is $G(\Lambda_s)$, the normalized second moment of a Voronoi cell of
$\Lambda_s$ ($=\frac{1}{12}$ for the
square lattice), $\nu_{2}=N_2$, $\nu_s=N_1N_2$ and ${N_{\cap}}/{[\Lambda_1:\Lambda_{\cap}]}=N_1$. Therefore 
\begin{eqnarray}
{J}_{s_1} &  = &  \frac{\gamma_1\gamma_2}{\gamma_1+\gamma_2}\frac{[\Lambda_1:\Lambda_{\cap}]}{N_{\cap}\nu_2} G(\Lambda_s)\nu_s^{1+2/L} \nonumber \\
\label{eq:Comp1}
 &  = &  \frac{\gamma_1\gamma_2}{\gamma_1+\gamma_2} G(\Lambda_s)\nu_s^{2/L} \\
\label{eq:lb_app_Lag}
 &  = &  \frac{\gamma_1\gamma_2}{\gamma_1+\gamma_2} G(\Lambda_s)(N_1N_2)^{2/L}.
\end{eqnarray}
If all the lattices in question are scaled by $\beta$, then 
\begin{equation}
\label{eq:lb_app_Lag:scaled}
{J}_{s_1} = \beta^2\frac{\gamma_1\gamma_2}{\gamma_1+\gamma_2}
G(\Lambda_s)(N_1N_2)^{2/L}.
\end{equation}

The rate of the $i$th description is given by 
\begin{equation}
\label{eq:rates}
R_i = h(p) - \frac{1}{L}\log_2(N_i) - \frac{1}{L}\log_2(\beta^L), \,\,\,\, i = 1,2.
\end{equation}
Therefore 
\begin{equation}
\label{eq:RatesVsIndex}
N_i = \frac{1}{\beta^L}2^{Lh(p)}2^{-LR_i}, \,\,\,\, i=1,2.
\end{equation}
The scale factor $\beta^2$ is related to the differential entropy of the source
and the rate $R_0$ through,
\begin{equation}
\label{eq:d0}
\beta^2=2^{2h(p)}2^{-2R_0}.
\end{equation}
Using (\ref{eq:RatesVsIndex},\ref{eq:d0}) in (\ref{eq:Comp1}) we obtain
\begin{eqnarray}
\label{eq:LagDist0}
{J}_{s_1} \approx \frac{\gamma_1\gamma_2}{\gamma_1+\gamma_2}
G(\Lambda_s)2^{2h(p)}2^{-2(R_1+R_2-R_0)}.
\end{eqnarray}
A bound for the term $J_{s_2}$ is obtained in terms of $\rho_{\cap}$, the covering radius of $\Lambda_{\cap}$,
by observing that for every $\lambda$, it is possible through a suitable $\Lambda_{\cap}$ shift,  
to satisfy
\begin{equation}
\left \|\lambda-\frac{\gamma_1\alpha_1(\lambda)+\gamma_2\alpha_2(\lambda)}{\gamma_1+\gamma_2}\right \|\leq \rho_{\cap}.
\end{equation}
Thus we have the inequality
\begin{equation}
J_{s_2}\leq (\gamma_1+\gamma_2)\rho^2_{\cap} \beta^2.
\label{eq:comp2}
\end{equation}
By comparing (\ref{eq:Comp1}) and (\ref{eq:comp2}) we observe that $J_{s_1}= \Theta(n^{4})$ whereas $J_{s_2}=\Theta(n^{2})$. Hence $J_{s_1}$ dominates $J_{s_2}$ and we obtain the approximation
\begin{eqnarray}
\label{eq:LagDist}
J \approx \frac{\gamma_1\gamma_2}{\gamma_1+\gamma_2}
G(\Lambda_s)2^{2h(p)}2^{-2(R_1+R_2-R_0)},
\end{eqnarray}
where $R_0$ determines the central distortion $\bar{d}_0$ and is given by
$\bar{d}_0=G(\Lambda)2^{2(h(p)-R_0)}$.

The approximations to the side distortions are obtained by using Figure \ref{fig:fig4} and
the following analysis. The channel 1 distortion is given by
\begin{equation}
\label{eq:SideDist1}
\displaystyle \bar{d}_1 = \frac{1}{N_{\cap}}
\sum_{\lambda\in V_{\Lambda_\cap}(0)}
||\alpha_1(\lambda)-\lambda||^2 + \bar{d}_0.
\end{equation}
Now the central approximation in the high rate analysis of the side distortions
is obtained by using Figure \ref{fig:fig4}.
The main idea being that the distance ${\rm AD}^2=||\lambda-\lambda_1||^2$ is well 
approximated by ${\rm AB}^2=||\bar{\lambda}-\lambda_1||^2$ at high rate.
We will formalize this notion below.
\begin{eqnarray}
\label{eq:SideDist2}
||\lambda-\lambda_1||^2 &=& {\rm AD}^2 = {\rm AC}^2 + {\rm CD}^2
\\ \nonumber
{\rm AF}^2 &\leq& {\rm AC}^2 \leq {\rm AG}^2
\\ \nonumber 
0 &\leq& {\rm CD}^2 \leq {\rm BD}^2
\end{eqnarray}

Now, by writing ${\rm BD}=\epsilon$ (as shown in Figure \ref{fig:fig4})
and using the geometry shown we obtain the following inequalities (note that
$\bar{\lambda}=(\gamma_1\lambda_1 +\gamma_2\lambda_2)/(\gamma_1+\gamma_2)$),

\begin{eqnarray}
\label{eq:SideDist3}
||\lambda_1-\bar{\lambda}||^2 \left [ 1 - \frac{\epsilon}{||\lambda_1-\bar{\lambda}||}
\right ]^2
\leq ||\lambda-\lambda_1||^2 \leq
||\lambda_1-\bar{\lambda}||^2 \left \{ \left [ 1 + \frac{\epsilon}
{||\lambda_1-\bar{\lambda}||} \right ]^2 + \left [ \frac{\epsilon}
{||\lambda_1-\bar{\lambda}||} \right ]^2 \right \}.
\end{eqnarray}

Therefore we can rewrite the above as
\begin{eqnarray}
\label{eq:SideDist4}
&||\lambda_1-\bar{\lambda}||^2  + ||\bar{\lambda}-\lambda||^2
-2||\lambda_1-\bar{\lambda}||||\bar{\lambda}-\lambda||
\leq ||\lambda-\lambda_1||^2 \leq \\ \nonumber
&||\lambda_1-\bar{\lambda}||^2 + ||\bar{\lambda}-\lambda||^2
+2||\lambda_1-\bar{\lambda}||||\bar{\lambda}-\lambda||
+||\bar{\lambda}-\lambda||^2
\end{eqnarray}
As $||\bar{\lambda}-\lambda||^2\geq 0$ and summing over $\lambda\in V_{\Lambda_{\cap}}(0)$,
we obtain
\begin{eqnarray}
\label{eq:SideDist5}
&\sum_{\lambda\in V_{\Lambda_{\cap}}(0)}\left [||\lambda_1-\bar{\lambda}||^2
-2||\lambda_1-\bar{\lambda}||||\bar{\lambda}-\lambda||\right ]
\leq \sum_{\lambda\in V_{\Lambda_{\cap}}(0)}||\lambda-\lambda_1||^2 \leq \\ \nonumber
&\sum_{\lambda\in V_{\Lambda_\cap}(0)}\left [||\lambda_1-\bar{\lambda}||^2 + 
2||\bar{\lambda}-\lambda||^2
+2||\lambda_1-\bar{\lambda}||||\bar{\lambda}-\lambda||\right ].
\end{eqnarray}
Now, by using the fact that $||\lambda_1-\bar{\lambda}||\leq \rho_{\cap}$,
we obtain
\begin{eqnarray}
\label{eq:SideDist6}
&\sum_{\lambda\in V_{\Lambda_{\cap}}(0)}||\lambda_1-\bar{\lambda}||^2\left[
1-2\frac{\sum_{\lambda\in V_{\Lambda_{\cap}}}||\lambda_1-\bar{\lambda}||
||\bar{\lambda}-\lambda||}
{\sum_{\lambda\in V_{\Lambda_{\cap}}}||\lambda_1-\bar{\lambda}||^2}\right ]
\leq \sum_{\lambda\in V_{\Lambda_{\cap}}(0)}||\lambda-\lambda_1||^2 \leq \\ \nonumber
&\sum_{\lambda\in V_{\Lambda_{\cap}}(0)}||\lambda_1-\bar{\lambda}||^2 \left [
1+2\frac{\sum_{\lambda\in V_{\Lambda_{\cap}}(0)}||\bar{\lambda}-\lambda||^2}
{\sum_{\lambda\in V_{\Lambda_{\cap}}(0)}||\lambda_1-\bar{\lambda}||^2}
+2\frac{\sum_{\lambda\in V_{\Lambda_{\cap}}(0)}||\lambda_1-\bar{\lambda}||
||\bar{\lambda}-\lambda||}
{\sum_{\lambda\in V_{\Lambda_{\cap}}(0)}||\lambda_1-\bar{\lambda}||^2}\right ].
\end{eqnarray}

Therefore, using these inequalities we obtain the following result.
\begin{lemma}
\label{lem:LCMapprox}
If $\gamma_1\neq 0, \gamma_2\neq 0$,
$\lim_{R_1\rightarrow \infty}\sum_{\lambda\in V_{\Lambda_{\cap}}(0)}||\lambda-\lambda_1||^2
=\sum_{\lambda\in V_{\Lambda_\cap}(0)}||\lambda_1-\bar{\lambda}||^2$
when $R_1-R_2 = C$ for some constant $C$.
\end{lemma}
\begin{proof}
For our sequence of lattices
\begin{equation}
\label{eq:norm_ineq}
0\leq \frac{\sum_{\lambda\in V_{\Lambda_{\cap}}(0)}||\lambda_1-\bar{\lambda}||
||\bar{\lambda}-\lambda||}
{\sum_{\lambda\in V_{\Lambda_{\cap}}(0)}||\lambda_1-\bar{\lambda}||^2} \leq P
\frac{\rho_{\cap}}{\nu_s^{1/L}},
\end{equation}
where $P$ is a constant that depends on $\Lambda_s$.
Hence, as $\frac{\rho_{\cap}}{\nu_s^{1/L}}=\Theta(n^{-1})$
we obtain $\lim_{R_1\rightarrow \infty}\frac{\rho_{\cap}}{\nu_s^{1/L}}=
\lim_{n\rightarrow \infty}\frac{\rho_{\cap}}{\nu_s^{1/L}}=0$---the desired result.
\end{proof}

Using Lemma \ref{lem:LCMapprox} we can write
\begin{equation}
\label{eq:SideDist7}
\sum_{\lambda\in V_{\Lambda_{\cap}}(0)}||\lambda-\lambda_1||^2 \approx 
\sum_{\lambda\in V_{\Lambda_{\cap}}(0)}||\lambda_1-\bar{\lambda}||^2 = 
\frac{\gamma_2^2}{(\gamma_1+\gamma_2)^2} \sum_{\lambda\in V_{\Lambda_{\cap}}(0)}
||\lambda_1-\lambda_2||^2
\end{equation}
at a high enough rate.
Therefore the side distortions are directly related to $J_s$ which was calculated
earlier. 
Hence the side distortions are approximated by
\begin{eqnarray}
\label{eq:SideDist}
\bar{d}_1 \approx \frac{\gamma_2^2}{(\gamma_1+\gamma_2)^2}G(\Lambda_s)2^{2h(p)}2^{-2(R_1+R_2-R_0)} 
\\ \nonumber
\bar{d}_2 \approx \frac{\gamma_1^2}{(\gamma_1+\gamma_2)^2}G(\Lambda_s)2^{2h(p)}2^{-2(R_1+R_2-R_0)}.
\end{eqnarray}
This allows us to find the approximate distortion ratio to be
$\frac{\bar{d}_1}{\bar{d}_2} \approx (\frac{\gamma_2}{\gamma_1})^2$ helping us to 
design the lattice quantizer
Though this approximation is asymptotic in the rate, we observed that it was
quite good for the quantizer design illustrated in the numerical results.

Next we examine the case when $\gamma_1=0,\gamma_2\neq 0$.
In this case we can show that,
\begin{eqnarray}
\label{eq:SideDist0}
\bar{d}_1 &\approx& G(\Lambda_s)2^{2h(p)}2^{-2R_1}
\\ \nonumber
\bar{d}_2 &\approx& G(\Lambda_s)2^{2h(p)}2^{-2(R_1+R_2-R_0)}.
\end{eqnarray}
Similarly the roles are reversed when $\gamma_1=0$ and $\gamma_2\neq 0$.

Let $R_0 = \frac{R_1+R_2}{2}(1+a)$, then 
$R_1+R_2-R_0=\frac{R_1+R_2}{2}(1-a)$.
Note that $a$ is chosen such that $a>\frac{|R_1-R_2|}{R_1+R_2}$ and therefore
$R_1+R_2-R_0<\min(R_1,R_2)$. 
Here we can clearly see the tradeoff in the central distortion 
$\bar{d}_0$ and the side distortions.

\subsection{Minimizing average distortion}
Suppose we know that the packet loss probability on channel 1 is
$p_1$ and the packet loss probability on channel 2 is $p_2$.
Then the average distortion is given by:
\begin{eqnarray}
\label{eq:Rsym_awgn}
\bar{D} = (1-p_1)(1-p_2)\bar{d}_0+(1-p_1)p_2 \bar{d}_1+(1-p_2)p_1 \bar{d}_2
+p_1p_2 \Expt[||\xbf||^2]
\end{eqnarray}
Now using the high rate approximations developed earlier,
we can find the optimal $\frac{\gamma_1}{\gamma_2}$
needed for minimizing the distortion.

\begin{claim}
The weights which minimize (\ref{eq:Rsym_awgn}) at high rate are given by
\begin{equation}
\label{eq:avgDist}
\frac{\gamma_1}{\gamma_2}=\frac{(1-p_1)p_2}{(1-p_2)p_1}.
\end{equation}
\end{claim}
{\rm 
\paragraph{Proof:}
To optimize (\ref{eq:Rsym_awgn}) we use the high rate expressions given in
(\ref{eq:SideDist}).
Using (\ref{eq:SideDist}) in (\ref{eq:Rsym_awgn}) we obtain
\begin{equation}
\label{eq:AvgDist1}
\bar{D} = A + B_1(\frac{\gamma_1}{\gamma_1+\gamma_2})^2 + 
B_2(\frac{\gamma_2}{\gamma_1+\gamma_2})^2,
\end{equation}
where $A,B_1,B_2$ do not depend on $\gamma_1,\gamma_2$ (they depend
on $R_1,R_2,R_0,\beta$).
Without loss of generality, we can use $\bar{\gamma}_1=\frac{\gamma_1}
{\gamma_1+\gamma_2}$ and $\bar{\gamma}_2=\frac{\gamma_2}{\gamma_1+\gamma_2}$.
Hence defining $\gamma=\bar{\gamma}_1=1-\bar{\gamma}_2$ and substituting
in (\ref{eq:AvgDist1}) we obtain
\begin{equation}
\label{eq:AvgDist2}
\bar{D} = A + B_1\gamma^2 + B_2(1-\gamma)^2.
\end{equation}
By differentiating (\ref{eq:AvgDist2}) with respect to $\gamma$
and setting it to zero we obtain the given result.
Note that this problem is convex (just differentiate (\ref{eq:AvgDist2}) 
twice and we see that it is always positive) and hence we have obtained
the minimum with respect to $\gamma$.
\hfill$\blacksquare$




\section{Numerical Results}
\label{sec:results}
In order to illustrate the performance of the proposed
quantizer, we evaluate its rate-distortion performance.
In order to compare its performance with that predicted by
information theory, we assume that there is an entropy (lossless)
coding of the quantizer output.
This is done for a Gaussian source with unit variance, for which
the multiple description rate-distortion problem
was solved by Ozarow \cite{OZAI}.

The example chosen is the $\ZZ^2$ lattice that we described in
Section \ref{sec:algo}.
The rates are chosen so that $R_1-R_2=\frac{1}{2}
\log_2(\frac{|\Lambda_2|}{|\Lambda_1|})$.

In Figure \ref{fig:d1-d2} we illustrate the tradeoff between
the two side distortions by varying $\gamma_1,\gamma_2$.

\begin{figure}[t]
\begin{center}
\includegraphics[scale=0.4]{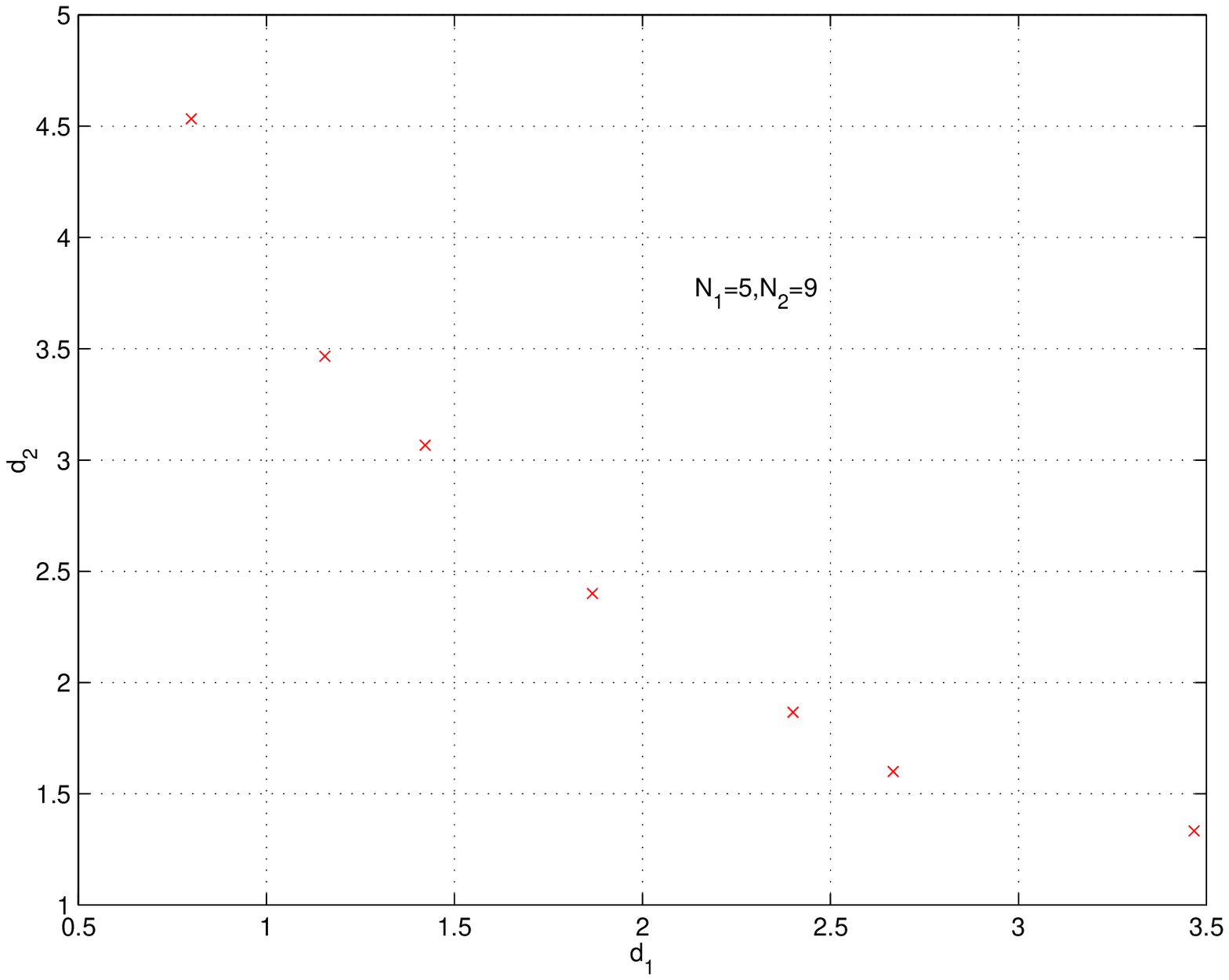}
\end{center}
\caption{Side distortions for fixed sublattices with varying 
$\gamma_1,\gamma_2$.}
\label{fig:d1-d2}
\end{figure}

In Figure \ref{fig:dR-gR} we have plotted a comparison
of $\frac{\gamma_2^2}{\gamma_1^2}$ with $\frac{\bar{d}_1}{\bar{d}_2}$.

\begin{figure}[t]
\begin{center}
\includegraphics[scale=0.4]{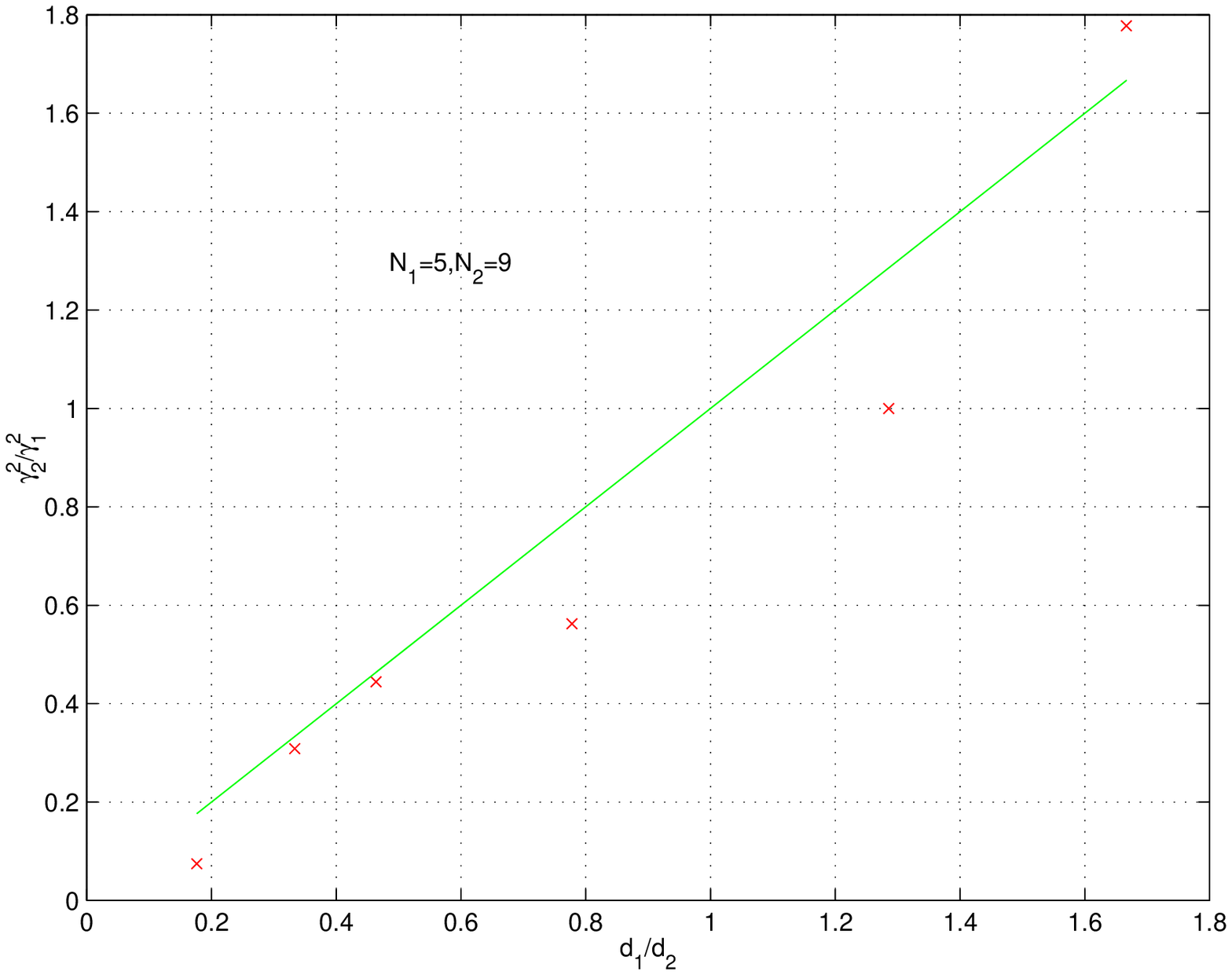}
\end{center}
\caption{}
\label{fig:dR-gR}
\end{figure}

In Figure \ref{fig:figRD} we have plotted the side distortions 
and compared them with those predicted by information theory \cite{OZAI}.
The key observation is that the distortion performance of the
lattice quantizer is approximately 3dB away from that
predicted by the rate-distortion bound.
This gap is due to the shaping gain that we will pick up when
we go to higher dimensions and using sublattices which have 
Voronoi cells which are close to spherical.
The $\ZZ^2$ lattice used in this example is more for illustrative purposes
and has very little shaping gain.

\begin{figure}[t]
\begin{center}
\includegraphics[scale=0.4]{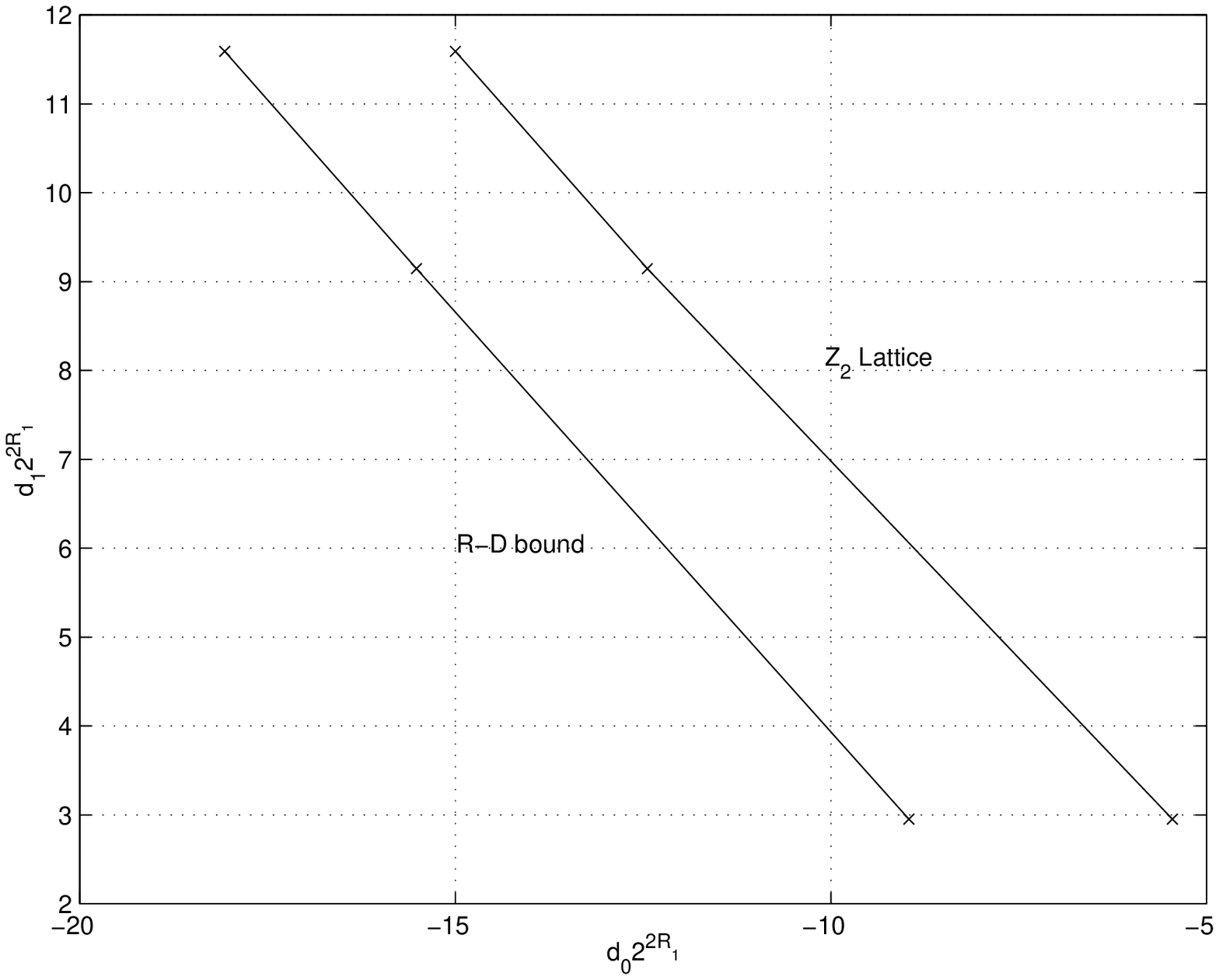}
\includegraphics[scale=0.4]{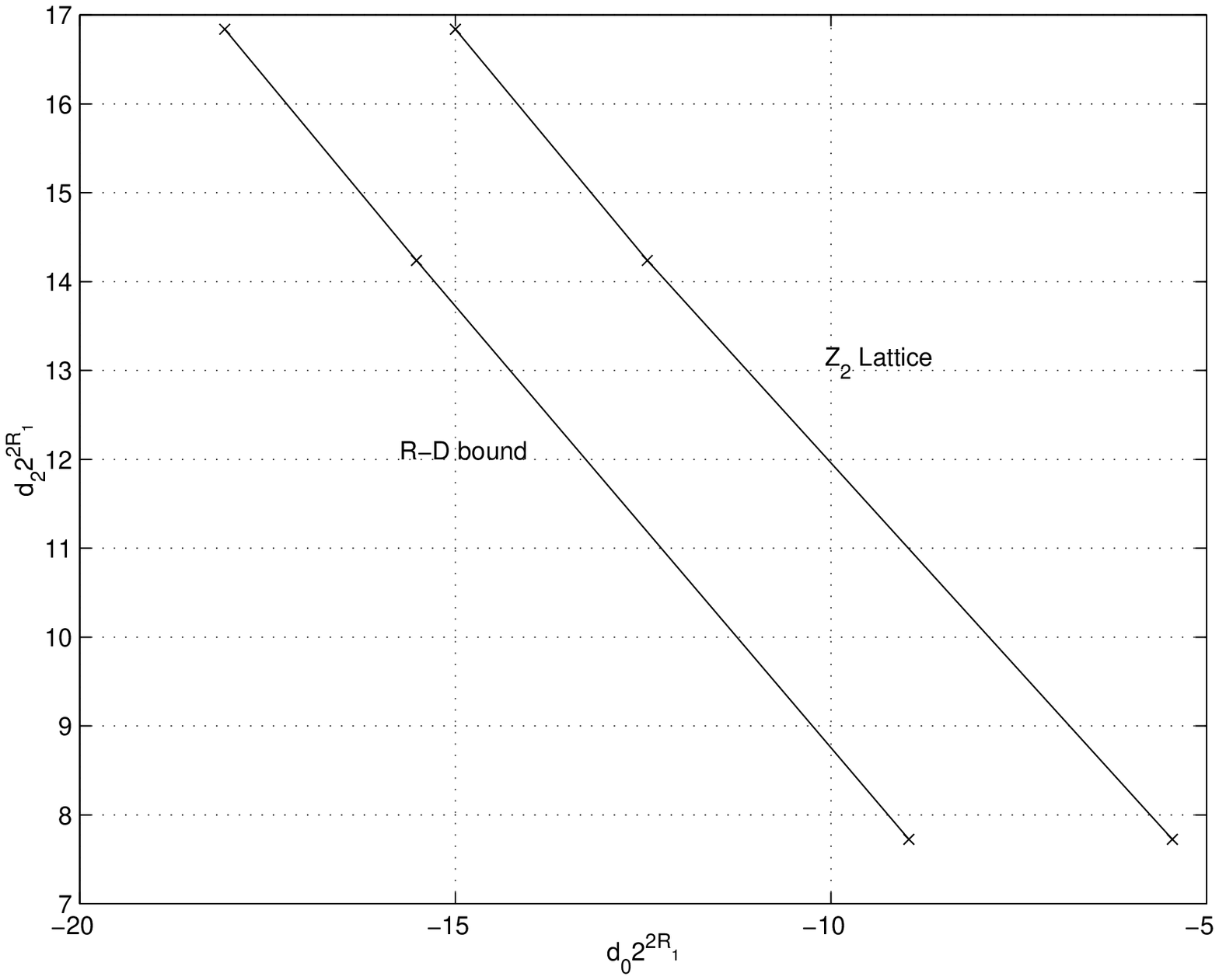}
\end{center}
\caption{Comparison of lattice distortion with rate-distortion bound.}
\label{fig:figRD}
\end{figure}

\section{Discussion}
\label{sec:conc}
In this paper we have designed asymmetric multiple description lattice quantizers.
This source coding scheme bridges the symmetric (balanced) multiple
description quantizers and completely hierarchical successive refinement
quantizers.
Though a lattice vector quantizer was illustrated, this scheme could also
be extended to other types of source coding schemes.


\small

\bibliography{/usr/suhas/bibs/bibrefs}
\bibliographystyle{ieeetr}

\appendix

\end{document}